\documentclass[10pt]{amsart}
\usepackage{amssymb, amsthm}\usepackage{amsmath}\usepackage{epsfig}
\usepackage{fancyhdr}

\theoremstyle{plain}
\newtheorem{Lem}{Lemma}[section]
\newtheorem{Not}[Lem]{Notation}
\newtheorem{Prop}[Lem]{Proposition}
\newtheorem{Cor}[Lem]{Corollary}

\newtheorem{The}[Lem]{Theorem}

\theoremstyle{definition}

\newtheorem{definition}[Lem]{Definition}
\newtheorem{Rem}[Lem]{Remark}
\newtheorem{Exe}[Lem]{Example}

\def\sn{W}
\def\renn{R(M)}
\def\crl{\Lambda}
\def\crlo{\Lambda_{\scriptscriptstyle \circ}}
\def\ZZe{\sn(e)}
\def\ZZf{\sn(f)}
\def\ZZse{\sn_\star(e)}

\newcommand{\GaTe}[1]{\textrm{Red}_\star(#1,\cdot)}
\newcommand{\Gae}[1]{\textrm{Red}(#1,\cdot)}
\newcommand{\DrTe}[1]{\textrm{Red}_\star(\cdot,#1)}
\newcommand{\Dre}[1]{\textrm{Red}(\cdot,#1)}

\def\ee{\varepsilon}
\def\C{\mathbb{C}}
\def\H{\mathcal{H}}
\begin{document}
\author{Godelle Eddy}
\address{Godelle Eddy\\ Universit\'e de Caen\\ Laboratoire Nicolas Oresme\\ 14032 Caen cedex France. email: eddy.godelle@math.unicaen.fr}
\date{\today}
\title{Generic Hecke Algebra for Renner monoids}
\maketitle
\begin{abstract} We associate with every Renner monoid~$R$ a \emph{generic Hecke algebra}~$\H(R)$ over~$\mathbb{Z}[q]$ which is a deformation of the monoid $\mathbb{Z}$-algebra of $R$. If $M$ is a finite reductive monoid with Borel subgroup~$B$ and associated Renner monoid~$R$, then we obtain the associated Iwahori-Hecke algebra~$\H(M,B)$ by specialising~$q$ in~$\H(R)$ and tensoring by~$\mathbb{C}$ over~$\mathbb{Z}$, as in the classical case of finite algebraic groups. This answers positively to a long-standing question of L.~Solomon. 

{\it 2000 Mathematics Subject Classification}: 20G40, 20C08, 20G05.\\keywords : Renner monoid, Hecke algebra.\end{abstract}
\section*{Introduction}
Consider the group~$\mathbb{G} = GL_n(\mathbb{F}_q)$ of invertible matrices over the finite field~$\mathbb{F}_q$. Denote by~$\mathbb{B}$ its subgroup of upper triangular matrices, and by~$\mathbb{T}$ its subgroup of diagonal matrices. Set~$\varepsilon = \frac{1}{|\mathbb{B}|}\sum_{b\in \mathbb{B}}b$ in~$\mathbb{C}[\mathbb{G}]$. The quotient group~$N_\mathbb{G}(\mathbb{T})/\mathbb{T}$ is isomorphic to the symmetric group~$S_n$. Moreover, the Iwohori-Hecke $\mathbb{C}$-algebra $\H(\mathbb{G},\mathbb{B}) = \varepsilon \mathbb{C}[\mathbb{G}]\varepsilon$ is isomorphic to $\oplus_{w\in S_n}\mathbb{C}w$ as a $\mathbb{C}$-vector space, and the structure constants in the multiplicative table lie in $\mathbb{Z}[q]$. More generally, if $G$ is a finite reductive group over~$\overline{\mathbb{F}}_q$, $B$ is a Borel subgroup of $G$, and $T$ a maximal torus included in~$B$, then~$N_G(T)/T$ is a Weyl group and the above results extend to the Hecke algebra $\H(G,B)$. Now, consider a \emph{finite reductive monoid}~$M$ over~$\overline{\mathbb{F}}_q$ as defined by Renner in~\cite{Ren2}. Such a monoid is a unit regular monoid and its unit group is a finite reductive group~$G$. Solomon introduced in~\cite{Sol} the notion of a \emph{Iwahori-Hecke algebra}~$\H(M,B)$ of a finite reductive monoid~$M$. Here,~$B$ is a Borel subgroup of~$G$. This $\mathbb{C}$-algebra is defined by $\H(M,B) = \varepsilon \mathbb{C}[M]\varepsilon$  where as before~$\varepsilon = \frac{1}{|B|}\sum_{b\in B}b$ in~$\mathbb{C}[M]$. In this framework, the Weyl group is replaced by an inverse monoid~$R$, which is called the \emph{Renner monoid} of~$M$. Its turns out that~$\H(M,B)$ is isomorphic to $\oplus_{r\in R}\mathbb{C}r$ as a $\mathbb{C}$-vector space. An isomorphism is given by $r\mapsto \tilde{T}_r = \sum_{x\in BrB}x$. Therefore, this is natural to address the question of the existence of a normalisation~$T_r = a_r\tilde{T}_r$ of the basis~$(\tilde{T}_r)_{r\in R}$ such that in this new basis~$(T_r)_{r\in R}$, the structure constants in the multiplicative table lie in $\mathbb{Z}[q]$ as in the case of finite reductive groups. Solomon considered this question in~\cite{Sol} and answered in the positive in the specific case where~$M = M_n(\mathbb{F}_q)$. In~\cite{Sol1}, he announced that in a forthcoming paper, he was going to extend his result and its proof to every finite reductive monoid that arises as the set of fixed points of a \emph{reductive monoid} over~$\overline{F}_q$ (see Section~\ref{sousect} for a definition) by the Frobenius map~$\sigma$ defined by $\sigma(x_{i,j}) = x_{i,j}^q$. But it seems that this result has never be published. In~\cite{Put} Putcha proves that for every finite reductive monoid, one can normalised the basis~$(\tilde{T}_r)_{r\in R}$ such that the structure constants become rational in~$q$. Howewer, the question remained open, and Renner concluded in~\cite[sec.~8.3]{Ren} that ``the delicate part here is obtaining integral structure constants''.  The main object of this article is to answer Solomon's question in the positive for every finite reductive monoid. We prove:    
\begin{The} 
Let $M$ be a finite reductive monoid over~$\overline{\mathbb{F}}_q$. Denote by $R$ the associated Renner monoid. There exists a normalisation of the basis~$(\tilde{T}_r)_{r\in R}$ of the Iwahori-Hecke algebra~$\H(M,B)$ such that the structure constants in the multiplicative table lie in~$\mathbb{Z}[q]$. Moreover, the coefficients of the polynomials depend on~$R$ only. \label{THintro}
\end{The}
In Section~\ref{sectionprincipale}, we provide explicit formulae (see Theorem~\ref{Th:genThe}), which are related to the existence of a length function on~$R$. Moreover, we deduce a finite presentation of~$\H(M,B)$ in the spirit of the classical presentation of~$\H(G,B)$ (see Corollary~\ref{THconclu} in Section~\ref{sousectionrappelamt}).

Mokler, Renner and Putcha consider families of monoids that are closed to reductive monoids (see \cite{Mok2,Mok3,Mok,Put2,Put3,PutRen} for instance. They are called \emph{finite monoids of Lie type} and  \emph{face monoids}. Indeed, finite reductive monoids are special cases of finite monoids of Lie type. To each of these groups can be associated a so-called \emph{Renner monoid}, whose properties are closed to Renner monoids of (finite) reductive monoids (See Examples~\ref{exe:ex2} and~\ref{exe:ex3} below).  This explains why these monoids are still called Renner monoids in the latter references. However, there is some differences between these monoids (see Remark~\ref{Rem:diffrenner} for a discussion). We introduce here the notion of a \emph{generalised Renner monoid}. All Renner monoids are examples of generalised Renner monoids. One motivation for this definition is to introduce a notion that plays for these various Renner monoids the role of the notion of a Coxeter system for Weyl group. We prove that all the properties shared by the various Renner monoids hold for generalised Renner monoid. In particular, it is a factorisable monoid and its unit group~$G$ is a Coxeter group. The crucial point regarding Solomon's question is that we can associate with each such generalised Renner monoid~$R$ a \emph{generic Hecke algebra}~$\H(R)$ which is a ring on the free $\mathbb{Z}[q]$-module with basis~$R$.  It turns out that Theorem~\ref{THintro} is a consequence of   
\begin{The} 
Let $M$ be a finite reductive monoid over~$\overline{\mathbb{F}}_q$ with Renner monoid~$R$. The Iwahori-Hecke algebra~$\H(M,B)$ is isomorphic to the $\mathbb{C}$-algebra~$\mathbb{C}\otimes_\mathbb{Z}\H_q(R)$, where $\H_q(R)$ is the specialisation of the generic Hecke algebra~$\H(R)$ at $q$. \label{THintro2}
\end{The}

The second main ingredient used in the proof of Theorem~\ref{THintro} is the existence of a length function~$\ell$ on every generalised Renner monoid~$R$. This length function is related to the canonical generating set~$S\cup\Lambda$, which equips every generalised Renner monoid. In the case of reductive monoids, we investigate the relation of this length function with the product of double classes. We prove in particular that  

\begin{Prop} Let $M$ be a reductive monoid with unit group~$G$ and Renner monoid~$R$. Fix a maximal torus~$T$ and a Borel subgroup~$B$ that contains~$T$ in $G$.\\(i) Let~$r$ lie in~$R$ and~$s$ lie in~$S$, then $$B s B r B = \left\{\begin{array}{ll}BrB,&\textrm{if } \ell(sr) = \ell(r);\\BsrB,&\textrm{if } \ell(sr) = \ell(r)+1;\\BsrB\cup BrB,&\textrm{if } \ell(sr) = \ell(r)-1.\end{array}\right.$$ \label{Pr:lienlongBBintro}
\\(ii) Let~$r$ lie in~$R$ and~$e$ lie in~$\Lambda$, then then $$B e B r B = BerB \textrm{ and } BrBeB = BreB$$ \end{Prop}

This result extends results obtained in~\cite{God,God2}, and leads to a similar result for finite reductive monoids. 

The paper is organised as follows. In Section~\ref{Sec:genheckalgr}, we introduce the notion of a \emph{generalised Renner monoid}, provide examples and investigate properties of such monoids. In particular, we define the length function~$\ell$ and prove that a \emph{generic Hecke algebra} can be associated with every generalised Renner monoid. In Section~\ref{sectionprincipale}, we first recall the notion of a reductive monoid and prove Proposition~\ref{Pr:lienlongBBintro}. Then we introduce the notion of a Iwahori-Hecke algebra in the context of monoid theory. We prove some motivating general results for such algebras. These results are probably well-known by semigroup experts, but we have not be able to find references for them. Finally, we turn to finite reductive monoids and conclude with the proof of Theorem~\ref{THintro} and~\ref{THintro2}.
\section{Generic Hecke algebra}
\label{Sec:genheckalgr}
The notion of a Coxeter group has been introduced in order to study Weyl groups. Our objective in this section is to develop a similar theory for Renner monoids. We need first to recall some standard notions and introduce useful notations.
\subsection{Basic notions and notations}
We refer to~\cite{How} for a general introduction on Semigroup Theory, and to~\cite{Fit} for a survey on factorisable inverse monoids.  We refer to~\cite{Bou} for general theory and proofs on Coxeter systems.

\subsubsection{Background on Semigroup Theory}
If $M$ is a monoid, we let $E(M)$ and $G(M)$  its idempotent set and its unit group. We see a (\emph{lower}) \emph{semi-lattice} as a commutative idempotent semigroup where $a\leq b$ iff $ab = ba = a$. In particular, $a\land b = ab$.  A semigroup is \emph{unit regular} if $M = E(M)G(M) = G(M)E(M)$, and it is \emph{factorisable} if it is unit regular and $E(M)$ is a semi-lattice. In this latter case $M$ is \emph{invertible}, that is for every $x$ in $M$ there exists a unique $y$ in $M$ such that $xyx= x$ (and therefore $yxy = y$).
\subsubsection{Background on Coxeter Group Theory}
\label{sousectionrappelcgtintro}
\begin{definition} Let $\Gamma$ be a finite simple labelled graph whose labels are positive integers greater or equal than~$3$. We let denote~$S$ the vertex set of $\Gamma$. We let $\mathcal{E}(\Gamma)$ denote the set of pairs $(\{s,t\},m)$ such that either $\{s,t\}$ is an edge of $\Gamma$ labelled by~$m$, or $\{s,t\}$ is not an edge of $\Gamma$ and $m=2$. When $(\{s,t\},m)$ belongs to~$\mathcal{E}(\Gamma)$, we let $|s,t\rangle^m$ denote the word~$sts\cdots$ of length~$m$. The \emph{Coxeter group}~$W(\Gamma)$ associated with~$\Gamma$ is defined by the following group presentation $$\left\langle S \left| \begin{array}{ll}s^2 = 1&s\in S\\ |s,t\rangle^m = |t,s\rangle^m &(\{s,t\},m)\in \mathcal{E}(\Gamma) \end{array}\right.\right\rangle$$
In this case, one says that the pair~$(W(\Gamma),S)$ is a \emph{Coxeter system}, and that $W$ is a Coxeter group. The Coxeter graph is uniquely defined by the Coxeter system. \end{definition}

\begin{definition} Let~$(W,S)$ be a \emph{Coxeter system}.\\ (i) Let $w$ belong to $\sn$. The \emph{length}~$\ell(w)$ of $w$ is the minimal integer $k$ such that $w$ has a word representative of length~$k$ on the alphabet~$S$. Such a word is called a \emph{minimal word representative} of~$w$.\\(ii) The subgroup~$W_I$ generated by a subset~$I$ of $S$ is called a \emph{standard parabolic subgroup} of $W$.\end{definition}
A key tool in what follows is the following classical result.
\begin{Prop}\cite{Bou} Let $(W,S)$ be a \emph{Coxeter system} with Coxeter graph~$\Gamma$.\\(i) For every~$I\subseteq S$, the pair $(W_I,I)$ is a Coxeter system. Its graph~$\Gamma_I$ is the full subgraph of $\Gamma$ spanned by~$I$. \\(ii)  For every~$I,J\subseteq S$ and every element $w\in W$ there exists a unique element~$\hat{w}$ of minimal length in the double-class $W_JwW_I$. Furthermore there exists $w_1$ in $W_I$ and $w_2$ in $W_J$ such that $w = w_2\hat{w}w_1$ with $\ell(w) = \ell(w_1)+\ell(\hat{w})+\ell(w_2)$.  \label{proppopa2} 
\end{Prop}
Note that~$(ii)$ holds, in particular, when $I$ or $J$ are empty. The element~$\hat{w}$ is said to be \emph{$(I,J)$-reduced}. In the sequel, we let $\textrm{Red}(I,J)$ denote the set of \emph{$(I,J)$-reduced} elements. Note also that the pair~$(w_1,w_2)$ is not unique in general, but it becomes unique if we require that $w_2\hat{w}$ is $(\emptyset,J)$-reduced (or that $\hat{w}w_1$ is $(I,\emptyset)$-reduced).

\subsection{Generalised Renner monoids}
\subsubsection{Generalised Renner-Coxeter System}
If $R$ is a factorisable monoid and $e$ belongs to $E(R)$ we let $W(e)$ and $W_\star(e)$ denote the subgroups defined by $$W(e) = \{w\in G(R)\mid we = ew\}$$ $$W_\star(e) = \{w\in G(R)\mid we = ew = e\}.$$ The unit group~$G(R)$ acts on $E(R)$ by conjugacy. 
\begin{definition} (i) An \emph{generalised Renner-Coxeter system} is a triple~$(R,\Lambda,S)$ such that
\begin{enumerate}
\item[(ECS1)] $R$ is a factorisable monoid;
\item[(ECS2)] $\Lambda$ is both a transversal of $E(R)$ for the action of $G(R)$ and a sub-semi-lattice;
\item[(ECS3)] $(G(R),S)$ is a Coxeter system;
\item[(ECS4)] for every pair $e_1\leq e_2$ in $E(R)$  there exists $w$ in $G(R)$ and $f_1\leq f_2$ in~$\Lambda$ such that $wf_iw^{-1} = e_i$ for~$i=1,2$;
\item[(ECS5)] for every $e$ in $\Lambda$, the subgroups~$W(e)$ and $W_\star(e)$ are standard Coxeter subgroups of $G(R)$;
\item[(ECS6)] the map $e\in\Lambda\mapsto \lambda^\star(e) = \{s\in S\mid se = es \neq e\}$ is not decreasing: $e\leq f \implies \lambda^\star(e)\subseteq \lambda^\star(f)$.
\end{enumerate}
In this case, we say that $R$ is a \emph{generalised Renner monoid}. Following the standard terminology for Renner monoids, we call the section~$\Lambda$ the \emph{cross section lattice} of $R$, and we define the~\emph{type map} of $R$ to be the map~$\lambda :\Lambda\to S$ defined by  $W(e) = W_{\lambda(e)}$.
\end{definition}

\begin{Not}
for $e$ in $\Lambda$, we set  $$\lambda_\star(e) = \{s\in S\mid se = es = e\}$$  $$W^\star(e) = W_{\lambda^\star(e)}$$  
\end{Not}

\begin{Rem} Assume $(R,\Lambda,S)$ is a generalised Renner-Coxeter system.\\  
(i) Since $W_\star(e)$ is a standard Coxeter subgroup of $W(e)$, we have $$W_\star(e) = W_{\lambda_\star(e)}.$$ Moreover, This is clear that $W_\star(e)$ is a normal subgroup of~$W(e)$. As a consequence, $$W(e) = W_\star(e)\times W^\star(e)\textrm{ and }\lambda(e) = \lambda_\star(e)\cup\lambda^\star(e).$$(ii) Below, several results can be proved without assuming Property~(ECS6). However this is a crucial tool in the proof of Theorem~\ref{Th:genThe} and Proposition~\ref{Prop:inverprese}.\\(iii) If $E(R)$ is finite and a lower semi-lattice, then it has to be a lattice. This is so for \emph{Renner monoids} associated with \emph{reductive monoids}.\\(iv) the map $\lambda_\star$ is not increasing: $$e\leq f \implies \lambda_\star(f)\subseteq \lambda_\star(e).$$(v) We can have $\lambda_\star(e) =  \lambda_\star(f)$ and $\lambda^\star(e) =  \lambda^\star(f)$ for $e\neq f$ (see \cite[Sec.~2.3]{God2}).     \label{remdepar}
\end{Rem}

Now we provide some examples of generalised Renner monoids.
\begin{Exe}  Let $M$ be a reductive monoid (see Section~\ref{sousect} for a definition, and Example~\ref{exerook}). The associated Renner\label{exe:ex1} monoid~$\renn$ of~$M$ is a generalised Renner monoid by \cite{Ren}. 
 \end{Exe} 
\begin{Exe}
Let $M$ be a \emph{abstract finite monoid of Lie type} (see~\cite{Put2},\cite{PutRen} or \cite{Ren} for a definition. Note that these groups are called \emph{regular split monoids} in~\cite{Put2}, and \emph{finite monoids of Lie type} in~\cite{PutRen}). The associated \emph{Renner monoid}~$\renn$ of~$M$ is a generalised Renner monoid. \label{exe:ex2} Property~(ECS6) follows from~\cite[Cor.~3.5(i)]{Put2}. The other defining properties hold by~\cite[Sec.~2]{Put3}. The seminal examples of an abstract finite monoid of Lie type is a Renner monoid of a \emph{finite reductive monoid}~\cite{Ren2}. In Section~3 we focus on these monoids.     
 \end{Exe}
\begin{Exe} \label{exe:ex3} Let~$G$ be a Kac-Moody group over a field~$\mathbb{F}$ of characteristic zero whose derived group is the special Kac-moody group introduced in~\cite{KaPe,KaPe2}. Denote by $(W,S)$ the associated Coxeter system. The Coxeter group~$W$ is infinite. Let~$Fa(X)$ be the set of \emph{faces} of its associated Tits cone~$X$ (see~\cite{Mok2} for details). The action of $W$ on $X$ induces an action on the lattice $Fa(X)$. The \emph{Renner monoid}~$R$ is the monoid~$W \ltimes Fa(X) /\sim$ where $\sim$ is the congruence on~$W\ltimes Fa(X)$ defined by $(w,R)\sim (w',R')$ if $R = R'$ and $w'^{-1}w$ fixes $R$ pointwise~\cite{Mok2}. Then $R$ is a generalised Renner monoid. Properties (ECS1), (ECS2), (ECS3) and (ECS5) are proved in~\cite{Mok2} (see also~\cite{Mok}). The  cross section lattice~$\Lambda$ can be identified with the set of infinite standard parabolic subgroups of $W$ that have no finite proper normal standard parabolic subgroups. The semi-lattice  structure is given by~$W_I\leq W_J$ if $J\subseteq I$. If $\Theta$ belongs to $\Lambda$, then $\lambda_\star(\Theta) = \Theta$ and $\lambda^\star(\Theta) = \{s\in S\mid \forall t\in \Theta, st = ts \}$. The latter equality clearly implies (ECS6). Finally, Property~(ECS4) can be deduced from~\cite[Theorem 2 and 4]{Mok}.  
 \end{Exe}
\begin{Rem} \label{Rem:diffrenner} In Examples~\ref{exe:ex1},~\ref{exe:ex2} and~\ref{exe:ex3} we provide examples of generalised Renner monoids that are all called \emph{Renner monoid} in the literature. From our point of view, this is not a suitable terminology since there is crucial differences between these monoids. Therefore, using the same terminology may be misleading. For instance, for Renner monoids of reductive monoids one has~$\lambda_\star(e) = \bigcap_{f\leq e} \lambda(f)$ and~$\lambda^\star(e) = \bigcap_{f\geq e} \lambda(f)$. This is not true in general for Renner monoids associated with abstract finite monoids of Lie type (see~\cite{PutRen} for a details). In Renner monoids of reductive monoids and of abstract finite monoids of Lie type, all maximal chains of idempotents have the same size. This is not true for Renner monoids of example~\ref{exe:ex3}, as explained in~\cite{Mok2}.
\end{Rem}

\subsubsection{Presentation for generalised Renner monoids}
For all this section, we fix a generalised Renner-Coxeter system~$(R,\Lambda,S)$. We let~$W$ denote the unit group of~$R$. Our objective is to prove that important properties shared by Renner monoids of Examples~\ref{exe:ex1}, \ref{exe:ex2}, \ref{exe:ex3} can be deduced from their generalised Renner-Coxeter system structure. In particular, we extend to this context the results obtained in~\cite{God2}. By Proposition~\ref{proppopa2}, For every~$w$ in~$\sn$ and every~$e,f$ in~$\Lambda$, each of the sets~$w\ZZe$, $\ZZe w$, $w\ZZse$, $\ZZse w$ and $\ZZe w \ZZf$ has a unique element of minimal length.  In order to simplify notation, we set $\Dre{e} = \textrm{Red}(\emptyset;\lambda(e))$, $\Gae{e} = \textrm{Red}(\lambda(e),\emptyset)$; $\DrTe{e} = \textrm{Red}(\emptyset,\lambda_\star(e))$; $\GaTe{e} = \textrm{Red}(\lambda_\star(e),\emptyset)$; $\textrm{Red}(e,f) = \textrm{Red}(\lambda(e),\lambda(f))$.

\begin{Prop}\label{fnrenner} For every~$r$ in~$R$,\\(i) there exists a unique triple~$(w_1,e,w_2)$ with~$e\in \crl$,~$w_1\in \DrTe{e}$ and~$w_2\in \Gae{e}$ such that~$r = w_1ew_2$;\\(ii) there exists a unique triple~$(v_1,e,v_2)$ with~$e\in \crl$,~$w_1\in \Dre{e}$ and~$w_2\in \GaTe{e}$ such that~$r = v_1ev_2$ \end{Prop}
Following~\cite{Ren}, we call the triple~$(w_1,e,w_2)$ the \emph{normal decomposition} of $r$.
\begin{proof} Let us prove~(i). The proof of~(ii) is similar.  Let $r$ belong to the monoid~$R$. By Property~(ECS1), there exists $e$ in~$E(R)$ and $w$ in $W$ such that $r = ew$. By Property~(ECS2) there exists $e_1$ in $\Lambda$ and $v$ in~$W$ such that $e = ve_1v^{-1}$. Then $r = vew_1$ with $w_1 = v^{-1}w$. By Remark~\ref{remdepar}(i), we can write $v = v_1v'_1$ and $w_1 = w'_2w''_2w_2$ with $v_1$, $w_2$, $v'_1$, $w'_2$ and $w''_2$ in $\DrTe{e}$, $\Gae{e}$, $W_\star(e)$, $W^\star(e)$ and $W_\star(e)$, respectively. Then we have $r = v_1w'_2ew_2$, and $v_1w'_2$ belongs to $\DrTe{e}$, still by Remark~\ref{remdepar}(i). Now assume~$r = w_1ew_2 = v_1fv_2$ with $e,f$ in $\Lambda$, $w_1,v_1$ in~$\DrTe{e}$ and~$\DrTe{f}$, respectively, and $w_2,v_2$ in~$\Gae{e}$ and in~$\Gae{f}$, respectively. Then $(w_1w_2)w_2^{-1}ew_2 = (v_1v_2)v_2^{-1}fv_2$. This implies $w_2^{-1}ew_2 = v_2^{-1}fv_2$ by~\cite{Fit}. As a consequence, $e = f$ and $v_2w_2^{-1}$ lies in $W(e)$. Since $v_2$ and $w_2$ both belong to $\Gae{e}$, we must have $v_2 = w_2$. Now, it follows that $w_1e = v_1e$ and $w_1^{-1}v_1$ lies in~$W_\star(e)$. This implies $w_1 = v_1$ in $\DrTe{e}$. 
 \end{proof}

\begin{Lem} Let $e,f$ belong to $\Lambda$ and $w$ lie in $\textrm{Red}(e,f)$.\\\label{lemclef}
(i) There exists $h$ in $\Lambda$ such that $w$ belongs to~$W(h)$ and $ewf = wh$.\\
(ii) The element~$w$ lies in~$W_\star(h)$. Therefore, $wh = h$.
\end{Lem}
Note that in the above lemma we have $h\leq e\wedge f = ef$. In the sequel the element~$h$ is denoted by~$e\wedge_wf$. 
\begin{proof} The proof is similar to~\cite[Prop~1.21]{God2}. (i) Consider the normal decomposition~$(w_1,h,w_2)$ of $ewf$. By definition~$w_1$ belongs to~$\DrTe{h}$ and $w_2$ belongs to~$\Gae{h}$. The element $w^{-1}ewf$ is equal to $w^{-1}w_1hw_2$ and belongs to $E(R)$. Since $w_2$ lies in~$\Gae{h}$, this implies that $w_3 = w_2w^{-1}w_1$ lies in $W_\star(h)$, and that $f\geq w_2^{-1}hw_2$. By Property~(ECS4), there exists $w_4$ in $\sn$ and $f_1,h_1$ in $\Lambda$, with $f_1\geq h_1$, such that $w_4^{-1}f_1w_4 = f$ and $w_4^{-1}h_1w_4 = w_2^{-1}hw_2$. Since $\Lambda$ is a cross section for the action of $\sn$, we have~$f_1 = f$ and~$h_1 = h$. In particular, $w_4$ belongs to~$W(f)$. Since $w_2$ belongs to $\Gae{h}$, we deduce that there exists $r$ in $W(h)$ such that $w_4 = rw_2$ with $\ell(w_4) = \ell(w_2)+\ell(r)$. Then $w_2$ lies in $W(f)$, too. Now, write $w_1 = w'_1w''_1$ where $w''_1$ lies in $W^\star(h)$ and $w'_1$ belongs to $\Dre{h}$. One has $ewf = w'_1hw''_1w_2$, and $w_1w''_2$ lies in~$\GaTe{h}$. By symmetry, we get that $w'_1$ belongs to $W(e)$. The element~${w'}^{-1}_1ww^{-1}_2$ is equal to $w''_1{w}_3^{-1}$ and belongs to $W(h)$. But, by hypothesis $w$ lies in~$\textrm{Red}(e,f)$. Then we must have $\ell(w''_1{w}_3^{-1}) = \ell({w'}^{-1}_1)+\ell(w)+\ell({w}^{-1}_2)$. Since $w''_1{w}_3^{-1}$ belongs to $W(h)$, it follows that $w'_1$ and $w_2$ belong to $W(h)$ too. This implies $w_2 = w'_1 = 1$ and $w = w''_1w_3^{-1}$. Therefore, $ewf = hw''_1 = hw = wh$.\\(ii) This is a direct consequence of the following fact: for $h,e$ in $\Lambda$ such that $h \leq e$, we have $W(h)\cap \Gae{e}\subseteq W_\star(h)$ and $W(h)\cap \Dre{e}\subseteq W_\star(h)$. Assume $w$ lies in~$W(h)\cap \Dre{e}$, then we can write $w = w_1w_2 = w_2w_1$ where $w_1$ lies in $W_\star(h)$ and $w_2$ lies in $W^\star(h)$. Since $h\leq e$, we have $\lambda^\star(h)\subseteq\lambda^\star(e)$ and $W^\star(h)\subseteq W^\star(e)$. Since $w$ belongs to~$\Dre{e}$, this implies~$w_2 = 1$. The proof of the second inclusion is similar.
 
\end{proof}

\begin{Cor} 
(i) For every chain $e_1\leq e_2\leq\cdots\leq e_m$ in $E(R)$  there exists $w$ in $G(R)$ and a chain $f_1\leq f_2\leq\cdots\leq f_m$ in~$\Lambda$ such that $wf_iw^{-1} = e_i$ for every index~$i$.\\
(ii) If $\Lambda$ has an infimum~$e$, then $\lambda(e) = S$.\\
(iii) For all $e,f$ in~$\Lambda$ and $w$ in~$\textrm{Red}(e,f)$, one has $$ewf  = \max\{h\in \Lambda\mid h\leq e,\ h\leq f,\ w\in W(h)\} = fw^{-1}e.$$
 \end{Cor}
In the case of Renner monoids of reductive monoids, the lattice~$\Lambda$ has an infimum~$e$ and $\lambda(e) = \lambda_\star(e) = S$. In other words, $e$ is a zero element of $R$.
\begin{proof}
 (i) Assume $w_1e_1w_1^{-1}\leq \cdots \leq w_me_mw_m^{-1}$. We prove the result by induction on $m$. For $m = 2$ this is true by Property (ECS4). Assume $m\geq 3$. By induction hypothesis, we can assume $w_2 = \cdots =w_m$. We can also also assume that $w_1$ belongs to~$\Dre{e_1}$. By hypothesis, we have $w_1e_1w_1^{-1}w_2e_2w_2^{-1} = w_1e_1w_1^{-1}$. We can write $w_1^{-1}w_2 = v_1v_3v^{-1}_2$ with $v_1$ in $W(e_1)$, $v_2$ in $W(e_2)$ and $v_3$ in $\textrm{Red}(e_1,e_2)$. Then $w_1e_1w_1^{-1}w_2e_2w_2^{-1} = w_1v_1e_1v_3e_2v^{-1}_2w_2^{-1}$. If $v_3\neq 1$, then we get a contradiction by Lemma~\ref{lemclef}(i) and Proposition~\ref{fnrenner}. Then $v_2 = 1$ and $e_1e_2=e_1$. It follows that  $w_1v_1 = w_2v_2$.  Write $v_1 = v_{1\star}v^\star_1$ and $v_2 = v_{2\star}v^\star_2$ with $v_{i\star}$ in $W_\star(e_i)$ and $v_i^\star$ in $W^\star(e_i)$. We have $w_1v_{1\star}v^{-1}_{2\star} = w_2v^\star_2{v_1^\star}^{-1}$. Since~$\lambda_\star(e_2)\subseteq \lambda_\star(e_1)$ and~$\lambda^\star(e_1)\subseteq \lambda^\star(e_2)$, we get  that~$v_{1\star}v^{-1}_{2\star}$ and~$v^\star_2{v^\star_1}^{-1}$ lie in $W(e_1)$ and $W(e_2)$, respectively. Then $w_1e_1w_1^{-1} = we_1w^{-1}$ and $w_2e_2w_2^{-1} = we_2w^{-1}$ with $w = w_1v_{1\star}v^{-1}_{2\star}$. But $W(e_2)\subseteq W(e_j)$ for $j\in\{2,\cdots,m\}$. Therefore, $w_2e_jw_2^{-1} = we_jw^{-1}$ for every~$j\geq 2$.\\ 
(ii) if $s\in S$ does not belong to~$\lambda(e)$, then $ese< e$ in $\Lambda$.\\
(iii) This is clear that~$e\wedge_wf$ lies in~$\{h\in \Lambda\mid h\leq e,\ h\leq f,\ w\in W(h)\}$. Now, if $h\in \Lambda$ verifies $h\leq e$, $h\leq f$, and $w\in W(h)$, then $h (e\wedge_wf) =  hw^{-1}(ewf) = w^{-1}hwf = hf = h$. Therefore, $h\leq ewf$. The last equality follows form the fact that $w^{-1}$ belongs to $\textrm{Red}(f,e)$. 
\end{proof}

\begin{Prop} For every $w$ in $W$, we fix an arbitrary reduced word representative~$\underline{w}$. We set $\crlo = \Lambda\setminus\{1\}$. The monoid~$R$ admits the monoid presentation whose generating set is~$S\cup \crlo$ and whose defining relations are:
\begin{center}\begin{tabular}{lll}
(COX1)&$s^2 = 1$,&$s\in S$;\\
(COX2)&$|s,t\rangle^m = |t,s\rangle^m$,&$(\{s,t\},m)\in \mathcal{E}(\Gamma)$;\\
(REN1)&$se = es$,&  $e\in \crlo$, $s\in \lambda^\star(e)$;\\
(REN2)&$se = es = e$,& $e\in \crlo$, $s\in \lambda_\star(e)$;\\
(REN3)&$e\underline{w}f = e\!\wedge_w\!f$,& $e,f\in \crlo$, $w\in \textrm{Red}(e,f)$. 
\end{tabular}\end{center} 
\label{proppres2}
\end{Prop} 
\begin{proof} This is clear that the relations stated in the proposition hold in $R$. Conversely, every element~$r$ in~$R$ has a unique representing word~$\underline{w}e\underline{v}$ such that $(w,e,v)$ is its normal decomposition, and this is immediate that every representing word of~$r$  on~$S\cup \crlo$  can be transformed into~$\underline{w}e\underline{v}$ using the given relations only. 
\end{proof}
\begin{Rem} (i) The above presentation is not minimal in general. Some of the relations of type~(REN3) can be removed (see the proof of~\cite[Theorem~0.1]{God2} and Remark~\ref{Rem:rem2} below).\\
(ii) The reader may verify that the result of Proposition~\ref{proppres2} and its proof still hold if we do not assume Property~(ECS6), except that Relation (REN3) must be replace by 
\begin{center}\begin{tabular}{lll}(REN3')&$e\underline{w}f = \underline{w} (e\wedge_wf)$,& $e,f\in \crlo$, $w\in \textrm{Red}(e,f)$.
\end{tabular}\end{center} Indeed, Lemma~\ref{lemclef}(i) still hold.
 \end{Rem}

One may wonder whether every monoid defined by a monoid presentation like in Proposition~\ref{proppres2}. The answer is positive under some necessary assumptions:

\begin{definition} A \emph{generalised Renner-Coxeter data} is $4$-uple~$(\Gamma,\crlo,\lambda_\star,\lambda^\star)$ such that~$\Gamma$ is a Coxeter graph with vertex set $S$, $\crlo$ is a lower semi-lattice and~$\lambda^\star$, $\lambda_\star$ are two maps from $\crlo$ to $S$ that verifies \begin{enumerate}
\item[(a)] for every $e$ in~$\crlo$, the graphs spanned by $\lambda_\star(e)$ and $\lambda^\star(e)$ in~$\Gamma$  are not connected, and $$e\leq f \Rightarrow \lambda_\star(f)\subseteq\lambda_\star(e)\textrm{ and }\lambda^\star(e)\subseteq\lambda^\star(f).$$
\item[(b)] for every $f,g$ in~$\crlo$ and every $w\in Red(f,g)$ the set $$\left\{e\in \crlo \mid e\leq f,\ e\leq g\textrm{ and }w\in W_{\lambda(e)}\right\}$$ has a greatest element, denoted by $f\wedge_w g$.
\end{enumerate} \label{Prop:inverpreseee}
with $\lambda(e) = \lambda_\star(e)\cup \lambda^\star(e)$ for $e\in\crlo$ and $\textrm{Red}(e,f) = \textrm{Red}(\lambda(e),\lambda(f))$ in the Coxeter group~$W(\Gamma)$ associated with~$\Gamma$.
\end{definition}

Note that properties (a) and (b) hold in every generalised Renner-Coxeter system. Actually, if~$\crlo$ is any lower semi-lattice such that all maximal chains are finite, then Assumption~(b) is necessarily verified. 

\begin{The} Assume~$M$ is a monoid. There exists a generalised Renner-Coxeter system~$(M,\Lambda,S)$ if and only if there exists a generalised Renner-Coxeter data~$(\Gamma,\crlo,\lambda_\star,\lambda^\star)$, where~$S$ is the vertex set of~$\Gamma$, such that $M$ admits the following monoid presentation 

\begin{center}\begin{tabular}{lll}
(COX1)&$s^2 = 1$,&$s\in S$;\\
(COX2)&$|s,t\rangle^m = |t,s\rangle^m$,&$(\{s,t\},m)\in \mathcal{E}(\Gamma)$;\\
(REN1)&$se = es$,&  $e\in \crlo$, $s\in \lambda^\star(e)$;\\
(REN2)&$se = es = e$,& $e\in \crlo$, $s\in \lambda_\star(e)$;\\
(REN3)&$e\underline{w}f = e\!\wedge_w\!f$,& $e,f\in \crlo$, $w\in \textrm{Red}(e,f)$. 
\end{tabular}\end{center} 
Where~$\underline{w}$  is an arbitrary fixed minimal representing word of~$w\in W(\Gamma)$.\\  In this case, $W(\Gamma)$ is canonically isomorphic to the unit group of $M$, and $\crlo$ embeds in $M$ with $\Lambda = \crlo\cup\{1\}$.
\label{Prop:inverprese}\end{The}

Note that given a generalised Renner-Coxeter data~$(\Gamma,\crlo,\lambda_\star,\lambda^\star)$, Relations (COX1) and (COX2) implies that the monoid~$M$ defined by the presentation stated in Theorem~\ref{Prop:inverprese} does not depend  on the chosen representing words~$\underline{w}$. Theorem~\ref{Prop:inverprese} follows from the following lemmas.

\begin{Lem} Consider a generalised Renner-Coxeter data~$(\Gamma,\crlo,\lambda_\star,\lambda^\star)$ and the monoid~$M$ defined by the presentation stated in Theorem~\ref{Prop:inverprese}. Then for every $f,g$ in~$\crlo$ and every $w\in \textrm{Red}(f,g)$,\label{lem:propreci}
\begin{enumerate}
\item[($b_1$)] $e\wedge_1 f = e\wedge f$ and $e\wedge_w f \leq e\wedge f$; 
\item[($b_2$)] $e\wedge_wf = f\wedge_{w^{-1}}e$; 
\item[($b_3$)]  $w\in W_{\lambda_\star(e\wedge_w f)}$.
\end{enumerate}  
\end{Lem}
\begin{proof} Properties~$(b_1)$ and $(b_2)$ are immediate consequences of Assumption~(b). Properties~$(b_3)$ follows from Assumption~(a). The main argument is like in the proof of Lemma~\ref{lemclef}(ii). If $w$ doesnot belong to $W_{\lambda_\star(e\wedge_w f)}$, then we can write $w = w_\star w^\star$ with $w_\star\in W_{\lambda_\star(e\wedge_w f)}$ and $w^\star\in W_{\lambda^\star(e\wedge_w f)}$. But $W_{\lambda^\star(e\wedge_w f)}\subseteq W_{\lambda^\star(f)}$ and $w$ lies in~$\textrm{Red}(e,f)$. Therefore, $w^\star = 1$.  
\end{proof}

\begin{Lem} Consider a generalised Renner-Coxeter data~$(\Gamma,\crlo,\lambda_\star,\lambda^\star)$ and the monoid~$M$ defined by the presentation stated in Theorem~\ref{Prop:inverprese}. Let $\textrm{F\!M\!}(S\cup \crlo)$ be the free monoid on~$S\cup \crlo$, and $\equiv$ be the congruence on $F\!M\!(S\cup \crlo)$ generated by the defining relations of~$M$. Hence by definition,~$M$ is equal to $\textrm{F\!M\!}(S\cup \crlo)/\!\equiv$.\\     
(i) If $\omega_1$ and $\omega_2$ are two words on~$S$ such that $\omega_1\equiv \omega_2$, then they represent the same element in $W(\Gamma)$.\\
(ii) If $e$ lie in $\crlo$ and $\omega$ lie in $\textrm{F\!M\!}(S\cup \crlo)$ with $e\equiv \omega$, then the word~$\omega$ is equal to~$\nu_1e_1\nu_2\cdots e_k\nu_{k+1}$ where for every~$i$ we have $e\leq e_i$ in $\crlo$ and $\nu_i$ are words on $S$ whose images in~$W(\Gamma)$ belong to $W_{\lambda(e)}$. Furthermore, the image of the word~$\nu_1\nu_2\cdots \nu_{k+1}$ in $W_{\lambda^\star(e)} = W_{\lambda(e)}/W_{\lambda_\star(e)}$ is trivial.\label{leminterm}
\end{Lem}
\begin{proof} In this proof we write $\omega_1 \dot{=} \omega_2$ if the two words $\omega_1$, $\omega_2$ are equals. If the words $\omega_1$ $\omega_2$  represent the elements $w_1, w_2$ in $M$, respectively, then $\omega_1 \dot{=} \omega_2$ implies $\omega_1\equiv \omega_2$ and $w_1 = w_2$. Conversely, $w_1 = w_2$ if and only if $\omega_1\equiv \omega_2$. Point~(i) is clear: if $\omega_1\equiv \omega_2$ then one can transform $\omega_1$ into $\omega_2$ using relations (COX1) and (COX2) only, since the words in both sides of Relations (REN1-3) contain letters in~$\crlo$. Let us prove~$(ii)$. Write $\omega_1\equiv_1 \omega_2$ if one can transform $\omega_1$ into $\omega_2$ by applying one defining relation of $M$ on $\omega_1$.  If $e\equiv \omega$, then there exists $\omega_0 \dot{=} e, \omega_1, \cdots, \omega_r \dot{=} \omega$ such that $\omega_i\equiv_1 \omega_{i+1}$. We prove the result by induction on $r$. If $r = 0$ we have nothing to prove. Assume $r\geq 1$. By induction hypothesis, $\omega_{r-1} \dot{=} \mu_1f_1\mu_2\cdots \mu_{j}f_j\mu_{j+1}$ with $e\leq f_i$ in $\crlo$ and $\mu_i$ is a word on $S$ whose image in~$W(\Gamma)$ belongs to $W_{\lambda(e)}$, and the image of the word~$\mu_1\mu_2\cdots \mu_{j+1}$ in $W_{\lambda^\star(e)} = W_{\lambda(e)}/W_{\lambda_\star(e)}$ is trivial. We deduce the result for $\omega \dot{=} \omega_r$ by considering case by case the type of the defining relation applied to $\omega_{r-1}$ to obtain~$\omega_r$. The cases where the relation is of one of the types (COX1), (COX2) or (REN1) are trivial. The case where the relation is of type (REN2) follows from Property~(a) in Definition~\ref{Prop:inverpreseee}: by induction hypothesis, one has $\lambda_\star(f_i)\subseteq \lambda_\star(e)\subseteq\lambda(e)$. Finally, the case where the relation is of type (REN3) follows from properties~(a) and (b) by Lemma~\ref{lem:propreci}. If the image~$u_i$ of~$\mu_i$ in $W(\Gamma)$ belongs to~$\textrm{Red}(f_{i-1},f_{i})$ with $\mu_i = \underline{u}_i$ and~$\omega \dot{=} \mu_1f_1\cdots \mu_{i-1} (f_{i}\!\wedge_{u_i}\!f_{i+1}) \mu_{i+1}f_{i+2}\cdots f_j\mu_{j+1}$ then $e\leq f_{i-1}\!\wedge_{u_i}\!f_i$. Conversely, if $\omega = \mu_1f_1\mu_2\cdots f_{i-1}\mu_ie_{i}\underline{u}_ie_{i+1}\mu_{i+1}\cdots \mu_{j}f_j\mu_{j+1}$ where $f_i = e_{i}\!\wedge_{u_i}\!e_{i+1}$ for $e_i,e_{i+1}$ in $\crlo$ and some $u_i$ in $\textrm{Red}(e_i,e_{i+1})$, then $e\leq f_i\leq e_i$ and $e\leq f_i\leq e_{i+1}$; Moreover,~$u_i$ belongs to $W_{\lambda_\star(f_i)}$, which is included in~$W_{\lambda_\star(e)}$. In all these cases the words~$\nu_1\nu_2\cdots \nu_{k+1}$ and~$\mu_1\mu_2\cdots \mu_{j+1}$ represent the same element in $W_{\lambda(e)}/W_{\lambda_\star(e)}$, which is trivial by induction hypothesis. 
\end{proof}

\begin{proof}[Proof of Theorem~\ref{Prop:inverprese}]
Consider a generalised Renner-Coxeter system~$(M,\Lambda,S)$. Denote by $\Gamma$ the Coxeter graph with vertex set~$S$ of the unit group of $M$, and set $\crlo = \Lambda\setminus \{1\}$. It follows from previous results that~$(\Gamma,\crlo,\lambda_\star,\lambda^\star)$ is a generalised Renner-Coxeter data, and by Proposition~\ref{proppres2} that $M$ has the required monoid presentation. Conversely, consider a generalised Renner-Coxeter data~$(\Gamma,\crlo,\lambda_\star,\lambda^\star)$ and let~$M$ denote the monoid defined by the presentation stated in Theorem~\ref{Prop:inverprese}. By Lemma~\ref{leminterm}(i), the subgroup of $M$ generated by $S$ can be identified with~$W(\Gamma)$. Lemma~\ref{leminterm}(ii) implies that $\crlo$ injects in~$M$, as a set. Let $e,f$ be in~$\crlo$. In $M$ one has  $ef = fe = e\wedge_1 f = e\wedge f$. Assume furthermore that $w$ lies in~$W(\Gamma)$. Lemma~\ref{leminterm}(ii) implies also that $(wew^{-1})f = wew^{-1}$ if and only if $e\leq f$ in~$\crlo$ and $w$ lie in $W_{\lambda(e)}$. Let $wew^{-1}$ and $vfv^{-1}$ be in $E(M)$ with $e,f$ in~$\crlo$. Write $w^{-1}v = v_1v_2v_3$ with $v_2$ in $\textrm{Red}(e,f)$, $v_1$ in $W_{\lambda(e)}$ and $v_3$ in $W_{\lambda(f)}$. Then $ev_2f = e\wedge_{v_2}f$ and $v_2$ lies in~$W_{\lambda_\star(e\wedge_{v_2}f)}$. We get, $$\displaylines{wew^{-1}vfv^{-1} = wv_1e\wedge_{v_2}fv_3v^{-1} = wv_1f\wedge_{v_2}^{-1}ev_3v^{-1} = wv_1v_2fv^{-1}_2ev_2v_3v^{-1} =\hfill\cr\hfill wv_1v_2v_3fv_3^{-1}v^{-1}_2v_1^{-1}ev_1v_2v_3v^{-1} = vfv^{-1}wew^{-1}.}$$ It is easy to see that every representing word~$\omega$ on $S\cup\crlo$ of en element~$w$ of $M$ can be transformed into a word $\omega_1e\omega_2 \equiv \omega_1e\omega_1^{-1}\omega_1\omega_2$ where $e$ belongs to~$\Lambda = \crlo\cup\{1\}$ and $\omega_1,\omega_2$ represent words in~$W(\Gamma)$. Moreover, if $\omega$ contains some letter in~$\crlo$, then~$e$ has to be in~$\crlo$. Therefore, $M$ is unit regular and $G(M) = W(\Gamma)$. In particular Property~(ECS3) holds. Assume $w = w_1ew_2$ lies in~$E(M)$ with $w_1,w_2$ in $W(\Gamma)$ and $e$ in~$\Lambda$. If $e = 1$ then $w_1w_2$ has to be equal to $1$ in~$W(\Gamma)$. Assume $e\neq 1$. Then $w_1ew_2w_1ew_2 = w_1ew_2$, and $ew_2w_1e = e$. By Lemma~\ref{leminterm}(ii), $w_2w_1$ belongs to $W_{\lambda_\star(e)}$ and~$w=w_1ew_1^{-1}$. Thus $E(M) = \{wew^{-1}\mid e\in\Lambda, w\in W(\Gamma)\}$ is a semi-lattice and Property~(ECS1) holds. Let $w_1,w_2,v_1,v_2$ be in~$W(\Gamma)$ and $e,f$ be in $\Lambda$  such that $w_1ew_2 =  v_1fv_2$ in $M$. Then $e =  w_1^{-1}v_1fv_2w_2^{-1}$ and $e\leq f$. By symmetry, $e = f$ and the elements~$w_1^{-1}v_1$ and~$v_2w_2^{-1}$ belong to~$W_{\lambda(e)}$. This implies that $\Lambda$ is a transversal of $E(M)$ for the action of $W(\Gamma)$ and a sub-semi-lattice of~$E(M)$. Therefore, we get Property~(ECS2). Furthermore, if $w_2 = v_1 = 1$ and $v_2 = w_1$, then $w_1$ lies in $W_{\lambda(e)}$. If $w_2 = v_1 = v_2 = 1$, then~$w_1$ lies in $W_{\lambda_\star(e)}$ by Lemma~\ref{leminterm}(ii). Property~(ECS5) follows. If $wew^{-1}\leq vfv^{-1}$, then $wew^{-1} vfv^{-1} = wew^{-1}$ and $ew^{-1}vfv^{-1}w= e$ Then $w^{-1}v$ lies in $W_{\lambda_\star(e)}\times W_{\lambda^\star(e)}$, which is included in $W_{\lambda_\star(e)}\times W_{\lambda^\star(f)}$. As a consequence, Property~(ECS4) holds. Finally, Property~(ECS6) holds by hypothesis. 
\end{proof}
\subsubsection{Length function for generalised Renner-Coxeter systems}
As explained in the introduction, to answer  Solomon's question, we need to define a length function on finite reductive monoids. Here we introduce this length function in the general context of generalised Renner-Coxeter systems. This extends results obtained in~\cite{God} and~\cite{God2}. As before,~$(R,\Lambda,S)$ is a generalised Renner-Coxeter system. The unit group of~$R$ is denoted by~$W$, and we set~$\crlo = \Lambda\setminus\{1\}$.\label{seclong}

\begin{definition}\label{deflenfun2} (i) We set~$\ell(s) = 1$ for~$s$ in~$S$ and~$\ell(e) = 0$ for~$e$ in~$\crl$. Let~$x_1,\ldots, x_k$ be in~$S\cup\crlo$  and consider the word~$\omega  =  x_1\cdots x_k$. Then the \emph{length} of the word~$\omega $ is the integer~$\ell(\omega )$ defined by~$\ell(\omega ) = \sum_{i = 1}^k\ell(x_i)$.\\
(ii) The \emph{length} of an element~$w$ which belongs to $R$ is the integer~$\ell(w)$ defined by $$\ell(w) =  \min\left\{\ell(\omega) \mid \omega \textrm{ is a word representative of }w \textrm{ over }S\cup\crlo\right\}.$$ If $\omega$ is a word representative of~$\omega$ such that~$\ell(w) = \ell(\omega)$, we say that~$\omega$ is a \emph{minimal word representative} of $w$.    \end{definition}
 
\begin{Prop} \label{propproplenght}Let~$r$ belong to~$R$.\\ (i) The length function~$\ell$ on~$R$ extends the length function~$\ell$ defined on~$\sn$.\\(ii)~$\ell(r) = 0$ iff~$r$ lies in~$\crl$.\\(iii) If~$s$ lies in~$S$ then~$|\ell(sr)-\ell(r)|\leq 1$.\\(iv) If~$r'$ belongs to~$R$, then~$\ell(rr')\leq \ell(r)+\ell(r')$.\end{Prop}
\begin{proof}This is direct consequences of the definition of the length function. \end{proof}

\begin{Prop}\label{propproplenght2}
Let~$r$ belong to~$R$. If $(w_1,e,w_2)$ is the normal decomposition of $r$, then $$\ell(r) = \ell(w_1)+\ell(w_2).$$
\end{Prop}
\begin{proof} Using the relations of the monoid presentation of~$R$ stated in Proposition~\ref{proppres2}, every representative word of $r$ can be transformed into~$\underline{w}_1e\underline{w}_2$ without increasing the length. Therefore~$\ell(r) = \ell(\underline{w}_1)+\ell(e)+\ell(\underline{w}_2) = \ell(w_1)+\ell(w_2)$.   \end{proof}

From the proof of the above proposition, we also deduce that
\begin{Cor}\label{cor:pasaugleng}
Let~$r$ belong to~$R$ and $\omega_1,\omega_2$ be two minimal word representatives of $r$. Using the relations of the monoid presentation of~$R$ stated in Proposition~\ref{proppres2}, one can transform~$\omega_1$ into~$\omega_2$ without increasing the length. 
\end{Cor}

\subsubsection{Matsumoto's Lemma for generalised Renner-Coxeter systems}
In this section we state and prove some technical results that play the role of \emph{Matsumoto's Lemma} in the context of generalised Renner-Coxeter systems. we need these results when proving Theorem~\ref{Th:genThe}. 
As before,~$(R,\Lambda,S)$ is a generalised Renner-Coxeter system. Let us first recall Matsumoto's Lemma.
\begin{Lem}\cite[Sec.~7.2]{Hum}\label{lemcoxbienconnu}
Consider a Coxeter system~$(W,S)$. Let $w$ belong to~$W$ and $s,t$ belong to $S$. If $\ell(swt) = \ell(w)$ and $\ell(sw) = \ell(wt)$, then $sw = wt$. 
\end{Lem}

\begin{Lem} \label{lem:ewwt1}Let $r$ belong to $R$ and $s,t$ belong to $S$. Let $(w_1,e,w_2)$ be the normal decomposition of $r$. Then\\(i) $\ell(sr) = \ell(r)\pm 1$ if and only if the normal decomposition of $sr$ is $(sw_1,e,w_2)$. In this case, $\ell(sr) - \ell(r) = \ell(sw_1) - \ell(w_1)$.\\(ii)  $\ell(sr) = \ell(r)$ if and only if $sr = r$ if and only if $sw_1 = w_1u$ for some $u$ in~$\lambda_\star(e)$. In this case, $\ell(sw_1) = \ell(w_1) +1$. \\(iii) $\ell(rt) = \ell(r)\pm 1$ if and only if the normal decomposition of $rt$ is either $(w_1,e,w_2t)$ or $(w_1u,e,w_2)$ for some $u$ in~$\lambda^\star(e)$. Furthermore, in the former case $\ell(rt) - \ell(r) = \ell(w_2t) -\ell(w_2)$, and in the latter case $w_2t = uw_2$ with $\ell(w_2t) = \ell(w_2)+1$.\\(iv)  $\ell(rt) = \ell(r)$ if and only if $r = rt$ if and only if $w_2t = uw_2$ for some $u$ in~$\lambda_\star(e)$.\\(v) If~$\ell(srt) = \ell(r)$ and $\ell(sr) = \ell(rt) \neq \ell(r)$, then there exists $u$ in~$\lambda^\star(e)$ such that $sw_1 = w_1u$ and $uw_2 = w_2t$. As a consequence, $sr = rt$.   \end{Lem}
\begin{proof} Recall that $|\ell (sr) -\ell (r)|\leq 1$ and $|\ell (rt) -\ell (r)|\leq 1$. The normal decomposition of $sr$ is $(sw_1,e,w_2)$ if and only if $sw_1$ belongs to $\DrTe{e}$. Since $w_1$ belongs to $\DrTe{e}$, this is clearly the case if $\ell(sw_1) = \ell(w_1)-1$. Assume~$\ell(sw_1) = \ell(w_1)+1$ and $sw_1$ does not belong to $\DrTe{e}$. Then we can write $sw_1 = w'_1u$ for some $u$ in $\lambda_\star(e)$ such that $\ell(sw_1) = \ell(w'_1)+1$. In particular, $\ell(sw_1u) = \ell(w'_1) = \ell(w_1)$. On the other hand, $\ell(w_1u) = \ell(w_1)+1 = \ell(sw_1)$ because $w_1$ belongs to $\DrTe{e}$, and $u$  lies in~$\lambda_\star(e)$. By Lemma~\ref{lemcoxbienconnu}, we get $sw_1 = w_1u$ and $sr = sw_1ew_2 = w_1uew_2 = w_1ew_2 = r$. This proves $(i)$ and $(ii)$ since the other implications are obvious. The normal decomposition of $rt$ is $(w_1,e,w_2t)$ if and only if $w_2t$ belongs to $\Gae{e}$. Since $w_2$ belongs to~$\Gae{e}$, this is clearly the case if $\ell(w_2t) = \ell(w_2)-1$. Assume~$\ell(w_2t) = \ell(w_2)+1$ and $w_2t$ does not belong to $\Gae{e}$. Then we can write $w_2t = uw'_2$ for some $u$ in $\lambda(e)$ such that $\ell(w_2t) = \ell(w'_2)+1$. As before we can conclude that $w_2t = uw_2$. If $u$ lies in $\lambda_\star(e)$ then $rt = r$. Otherwise, $u$ belongs to $\lambda^\star(e)$ and $w_1u$ belongs to $\DrTe{e}$. This is true since $u$ belongs to~$\lambda^\star(e)$ and therefore commutes with each element of~$\lambda_\star(e)$. Then the normal decomposition of~$rt$ is $(w_1u,e,w_2)$. This proves $(iii)$ and $(iv)$. Now assume $\ell(srt) = \ell(r)$ and $\ell(sr) = \ell(rt) \neq \ell(r)$. We claim that $\ell(w_2t) = \ell(w_2)+1$ and there exists $u$ in $\lambda(e)$ such that $uw_2 = w_2t$. If it was not the case, by above arguments, the normal decomposition of $srt$ would be $(sw_1,e,w_2t)$ and $\ell(srt) = \ell(r) \pm 2$. Since we assume~$\ell(rt) \neq \ell(r)$, the element~$u$ has to belong to~$\lambda^\star(e)$. Finally, using that $\ell(sr) = \ell(rt)\neq \ell(r) = \ell(srt)$ we deduce that $\ell(sw_1) = \ell(w_1u)$ and $\ell(w_1) = \ell(sw_1u)$, which in turn implies~$sw_1 = w_1u$ by Lemma~\ref{lemcoxbienconnu}.  
\end{proof}

\begin{Lem}\label{lem:ewwt2}Let $r$ belong to $R$,~$s$ belong to~$S$ and $f$ belong to $\crl$. Let $(w_1,e,w_2)$ be the normal decomposition of $r$.\\(i) If $\ell(rf) = \ell (r)$ then $w_2$ belongs to~$W(f)$.\\(ii)  If $\ell(fr) = \ell (r)$ then $w_1 = w'_1w''_1$ where $w'_1$ lies in ~$W(f)$ and $w''_1$ lies in $W^\star(e)$.\\(iii) If $\ell(sr) = \ell(r)-1$, then $\ell(srf)\leq\ell(rf)$. If $\ell(sr) = \ell(r)+1$, then $\ell(srf) \geq  \ell(rf)$. \\(iv) If $\ell(rs) = \ell(r)-1$, then $\ell(frs)\leq\ell(fr)$. If $\ell(rs) = \ell(r)+1$, then $\ell(frs) \geq  \ell(fr)$.
\end{Lem}
\begin{proof}  By definition of the normal decomposition, $w_2$ belongs to~$\Gae{e}$. Write $w_2 = w'_2w''_2$ with $w'_2,w''_2$ in the unit group~$W$ of~$R$ such that $\ell(w_2) = \ell(w'_2)+\ell(w''_2)$, $w''_2$ belongs to~$W(f)$ and $w'_2$ belongs to~$\Dre{f}$.  Then $w'_2$ lies in $\textrm{Red}(e,f)$. By Relation~(REN3), we  have~$rf = w_1(e\wedge_{w'_2}f)w''_2$. It follows that $\ell(w'_2) = 0$, and $w_2 = w''_2$. This proves~(i). The prove of (ii) is similar except that we need first to decompose $w_1$ in $w'_1w''_1$ where $w''_1$ lies in $W^\star(e)$ and $w'_1$ lies in~$\Dre{e}$.\\
(iii) Assume $\ell(sr) = \ell(r)-1$. Write $w_1 = sv_1$ with $\ell(w_1) = \ell(v_1)+1$, and write $w_2 = w'_2w''_2v'''_2$ with $w'_2,w''_2,w'''_2$ in $W$ such that $\ell(w_2) = \ell(w'_2)+\ell(w''_2)+\ell(w'''_2)$, where $w''_2$ belongs to~$W^\star(f)$, $w'''_2$ belongs to~$W_\star(f)$ and $w'_2$ belongs to~$\textrm{Red}(e,f)$. Then $(v_1,e,w_2)$ is the normal decomposition of~$sr$. One has $srf = v_1ew'_2fw''_2 = v_1e'w''_2$ where $e' = e\land_{w'_2}f$ belongs to~$\crl$. Write $w''_2 = v''_2v'_2v_2$ such that $\ell(w''_2) = \ell(v''_2)+\ell(v'_2)+\ell(v_2)$ with $v''_2\in W_\star(e')$, $v'_2\in W^\star(e')$ and $v_2\in \Gae{e'}$. We claim that $v''_2 = 1$. Indeed $w'_2$  belongs to~$W_\star(e')$ by Lemma~\ref{lemclef}(ii), and $w_2 = w'_2v''_2v'_2v_2w'''_2 = v'_2w'_2v''_2v_2w'''_2$ with $\ell(w_2) = \ell(v_2)+\ell(v'_2)+\ell(w'_2)+\ell(v''_2)+\ell(w'''_2)$. But $v'_2\in W^\star(e')\subseteq W^\star(e)$, since~$e'\leq e$ by Property~(ECS6), whereas $w_2$ belongs to~$\Gae{e}$ by definition of the normal decomposition. Hence, $v'_2 = 1$. Now, write~$v_1 = v'_1v''_1$ such that $\ell(v_1) = \ell(v'_1)+\ell(v''_1)$ with $v'_1\in \DrTe{e'}$ and $v''_1\in W_\star(e')$. Then $srf = v'_1e'v_2$ and $(v'_1,e',v_2)$ is the normal decomposition of $swf$. Since $\ell(ssr) = \ell(sr)+1$, we have $\ell(sv'_1v''_1) = \ell(sv_1) = \ell(v_1)+1$ by Lemma~\ref{lem:ewwt1}(i). This implies $\ell(sv'_1) = \ell(v'_1)+1$ and we cannot have $\ell(ssrf) =  \ell(srf)-1$, still by Lemma~\ref{lem:ewwt1}(i). Assume $\ell(sr) = \ell(r)+1$. let $(v_1,e,w_2)$ be the normal decomposition of~$r$, and $(v'_1,e',v_2)$ be the normal decomposition of~$rf$. It follows from above arguments that $v'_1$ left divides $v_1$. We conclude using Lemma~\ref{lem:ewwt1}: $\ell(sr) = \ell(r)+1\Rightarrow  \ell(sv_1) = \ell(v_1)+1 \Rightarrow  \ell(sv'_1) = \ell(v'_1)+1 \Rightarrow  \ell(srf) \geq  \ell(rf)$.  The proof of (iv) is similar.    

\end{proof}

\subsection{Free module over $R$}
For all this section, we assume~$(R,\Lambda,S)$ is a generalised Renner-Coxeter system. We let $W$ denote the unit group of $R$, and set $\crlo = \Lambda\setminus\{1\}$. We fix an arbitrary unitary associative ring~$A$. We let~$V$ denote the free $A$-module with basis elements $T_r$ for~$r\in R$.

\begin{The} Fix $q$ in $A$. There exists a unique structure of unitary associative $A$-algebra on~$V$ such that~$T_1$ is the unity element and the following conditions hold for every $x$ in $S\cup\crlo$ and every $r$ in $R$:\label{Th:genThe}
\begin{center}\begin{tabular}{ll}$T_xT_r = T_{xr}$,&if $x\in S$ and $\ell(xr) =  \ell(r) +1 $;\\
$T_xT_r = qT_{r}$,&if $x\in S$ and $\ell(xr) = \ell(r)$;\\$T_xT_r = (q-1)T_r+qT_{xr}$,&if $x\in S$ and $\ell(xr) = \ell(r) - 1;$\\ $T_xT_r = q^{\ell(r)-\ell(xr)}T_{xr}$,&if $x\in \crlo$.
\\\end{tabular}\end{center} 
\end{The}

 We follow the method explained in~\cite[Sec.~7.1]{Hum1} for the Hecke algebra of Coxeter groups. Let~$\mathcal{E} = End_A(V)$ the $A$-algebra of endomorphisms of the $A$-module~$V$. For~$s$ in~$S$ and $r$ in~$R$, we define $\rho_s$ in~$\mathcal{E}$ by \begin{center}\begin{tabular}{ll}$\rho_s(T_r) = T_{sr}$, &if $\ell(sr) =  \ell(r)+1$;\\$\rho_s(T_{r}) = qT_{r}$, &if $\ell(sr) =  \ell(r)$;\\$\rho_s(T_{r}) = (q-1)T_r+qT_{sr}$, &if $\ell(sr)  = \ell(r) -1$.\end{tabular}\end{center} For~$e$ in~$\crl$ and $r$ in~$R$, we define $\rho_e$ by $$\rho_e(T_r) =  q^{\ell(r)-\ell(er)}T_{er}$$ Similarly, for~$s$ in~$S$ and $r$ in~$R$, we define~$\overline{\rho}_s$ in~$\mathcal{E}$ by \begin{center}\begin{tabular}{ll}$\overline{\rho}_s(T_r) = T_{rs}$,&if $\ell(sr) =  \ell(r)+1$;\\$\overline{\rho}_s(T_{r}) = qT_{r}$,&if $\ell(r) =  \ell(rs)$;\\$\overline{\rho}_s(T_{r}) = (q-1)T_r+qT_{rs}$,&if $\ell(sr)  = \ell(r)-1$.\end{tabular}\end{center} For~$e$ in~$\crl$ and $r$ in~$R$, we define $\overline{\rho}_e$ by $$\overline{\rho}_e(T_r) =  q^{\ell(r)-\ell(re)}T_{re}.$$  
The key tool in the proof of Theorem~\ref{Th:genThe} is the following result.
\begin{Lem}\label{lem:commut} For every $x,y$ in $S\cup \crl$, $$\rho_x\overline{\rho}_y = \overline{\rho}_y\rho_x.$$ 
\end{Lem}

\begin{proof} Let $r$ belong to~$R$ and $x,y$ belong to $S\cup\crl$. We prove that~$\rho_x(\overline{\rho}_y (T_r))= \overline{\rho}_y(\rho_x(T_r))$. Clearly we can assume $x\neq 1$ and $y\neq 1$. By Proposition~\ref{propproplenght}, $\ell(xry)\leq \ell(x)+\ell(r)+\ell(y)\leq(r)+2$. We provide case by case as in~\cite{Hum}.\\
\underline{Case $1$: $\ell(xry) = \ell(r)+\ell(x)+\ell(y)$}. We must have $\ell(xr) = \ell(r)+\ell(x)$, $\ell(ry) = \ell(r)+\ell(y)$ and $\ell(xry) = \ell(ry)+\ell(x) = \ell(xr)+\ell(y)$. Therefore~$\rho_x(\overline{\rho}_y (T_r)) = \rho_x(T_{ry})= T_{xry} = \overline{\rho}_y(T_{xr}) = \overline{\rho}_y(\rho_x(T_r))$.\\

\noindent\underline{Case~$2$: $\ell(xry) = \ell(r)+1$}. We must have~$\ell(xr)\geq\ell(r)$,~$\ell(ry)\geq\ell(r)$, and~$x$ or~$y$, possibly both, belongs to $S$. If~$x$ or $y$ belongs to~$\crlo$, we are in Case~$1$. So we assume $x$ and $y$ belong to~$S$.\\
{\it Subcase~1:} $\ell(xr) = \ell(r)$, that is $xr = r$. Then $\ell(ry) = \ell(xry) = \ell(r)+1$ and $\ell(xry) = \ell(xr)+1$. Therefore~$\rho_x(\overline{\rho}_y (T_r)) = \rho_x(T_{ry})= qT_{xry} = \overline{\rho}_y(qT_{xr}) = \overline{\rho}_y(\rho_x(T_r))$. The case~$\ell(ry) = \ell(r)$ is similar. \\
{\it Subcase~2:}~$\ell(ry) = \ell(xr)= \ell(r)+1$. Then~$\ell(ry) = \ell(xr)= \ell(xry)$. We deduce that~$\rho_x(\overline{\rho}_y (T_r)) = \rho_x(T_{ry})= qT_{xry} = \overline{\rho}_y(T_{xr}) = \overline{\rho}_y(\rho_x(T_r))$.\\

\noindent\underline{Case~$3$: $\ell(xry) = \ell(r)$}.   If~$x$ and $y$ belong to~$\crlo$, we are in Case~$1$. So we assume this  is not the case.\\{\it Subcase~1:}~$x$ and $y$ belong to~$S$. Consider first the case $\ell(xr) = \ell(r)$. Then $xr = r$ and $\ell(xry) = \ell(ry) = \ell(r)$. Therefore,~$\rho_x(\overline{\rho}_y (T_r)) = \overline{\rho}_y(\rho_x(T_r)) = q^2T_r$. Assume now~$\ell(xr) \neq \ell(r)$. This implies $\ell(ry)\neq \ell(y)$ by symmetry. If $\ell(xr) = \ell(ry)$, by Lemma~\ref{lem:ewwt1}$(v)$ we have $xr = ry$. Hence, if $\ell(xr) = \ell(ry) = \ell(r) +1$, we have~$\overline{\rho}_y(\rho_x(T_r)) = \overline{\rho}_y(T_{xr}) = (q-1)T_{xr}+qT_{xry}$ and~$\rho_x(\overline{\rho}_y (T_r)) = \rho_x(T_{ry}) = (q-1)T_{ry}+qT_{xry}$. If $\ell(xr) = \ell(ry) = \ell(r) - 1$, we have~$\overline{\rho}_y(\rho_x(T_r)) = \overline{\rho}_y((q-1)T_{r}+qT_{xr}) = (q-1)T_{yr}+qT_{xry}$ and~$\rho_x(\overline{\rho}_y (T_r)) = \rho_x((q-1)T_{r}+qT_{ry}) = (q-1)T_{xr}+qT_{xry}$. Consider now the case  $\ell(xr) = \ell(r)+1$ and $\ell(ry) = \ell(r)-1$. Then~$\overline{\rho}_y(\rho_x(T_r)) = \overline{\rho}_y(T_{xr}) = (q-1)T_{xr}+qT_{xry} = \rho_x((q-1)T_{r}+q T_{ry}) = \rho_x(\overline{\rho}_y (T_r))$. The case where $\ell(xr) = \ell(r)-1$ and $\ell(ry) = \ell(r)+1$ is similar. \\{\it Subcase~2:}~$x$ belongs to~$S$ and~$y$ belong to~$\crlo$. We must have $\ell(xr)\geq \ell(r)$. Assume first~$\ell(xr) = \ell(r)$. We have $xr = r$ and $\ell(xry) = \ell(ry) = \ell(r)$. We get, $\overline{\rho}_y(\rho_x(T_r)) = \overline{\rho}_y(qT_r)) = q^{1+\ell(r)-\ell(ry)}T_{ry} =q^{\ell(r)-\ell(ry)}\rho_x(T_{ry}) = \rho_x(\overline{\rho}_y (T_r))$. Assume now~$\ell(xr) = \ell(r)+1$, then~$\overline{\rho}_y(\rho_x(T_r)) = \overline{\rho}_y(T_{xr}) =  q^{\ell(xr)-\ell(xry)}T_{xry} = q T_{xry}$. If $\ell(ry) = \ell(r)$ then $\ell(xry) = \ell(ry)$ and~$\rho_x(\overline{\rho}_y (T_r)) = \rho_x(T_{ry}) = qT_{xry}$. If $\ell(ry) < \ell(r)$, then $\ell(xry) = \ell(r) = \ell(ry) +1$ and ~$\rho_x(\overline{\rho}_y (T_r)) = q\rho_x(T_{ry}) = qT_{xry}$. The case $x\in \crlo$ and $y\in S$ is similar.\\

\noindent\underline{Case~$4$: $\ell(xry) < \ell(r)$}.\\
{\it Subcase~1:}~$x,y$ belong to~$\crlo$. Clearly,  $\rho_x(\overline{\rho}_y (T_r)) = \overline{\rho}_y(\rho_x(T_r)) = q^{\ell(r)-\ell(xry)}T_{xry}$.\\   
{\it Subcase~2:}~$x$ belongs to~$S$, $y$ belongs to~$\crlo$ and $\ell(xr) = \ell(r)$. Then $xr = r$ and $xry = ry$. This case is similar to the first case in Case~$3$ Subcase~$2$.\\   
{\it Subcase~3:}~$x$ belongs to~$S$, $y$ belongs to~$\crlo$ and~$\ell(xr) = \ell(r) - 1$. Applying Lemma~\ref{lem:ewwt2}, we get~$\ell(xry) \leq \ell(ry)$. We have $\overline{\rho}_y(\rho_x(T_r)) =$ $\overline{\rho}_y((q-1)T_r+qT_{xr}) =$ $(q-1)q^{\ell(r)-\ell(ry)}T_{ry}+q^{1+\ell(xr)-\ell(xry)}T_{xry}$ and $(\overline{\rho}_y (T_r)) = q^{\ell(r)-\ell(ry)} \rho_x(T_{ry})$.\\ Assume first $\ell(xry) = \ell(ry)-1$. Then $\ell(xr)-\ell(xry) = \ell(r)-\ell(ry)$ and $(\overline{\rho}_y (T_r)) = (q-1)q^{\ell(r)-\ell(ry)} T_{ry} + q^{1+\ell(r)-\ell(ry)} T_{xry}$.\\ Assume secondly that $\ell(xry) = \ell(ry)$, that is $xry = ry$. In this case, $(\overline{\rho}_y (T_r)) = q^{1+\ell(r)-\ell(ry)}T_{xry}$. But $1+\ell(xr)-\ell(xry) = \ell(r)-\ell(ry)$, therefore $\overline{\rho}_y(\rho_x(T_r)) = q^{1+\ell(r)-\ell(ry)} T_{ry}$.\\
{\it Subcase~4:}~$x$ belongs to~$S$, $y$ belongs to~$\crlo$ and $\ell(xr) = \ell(r)+1$. By Lemma~\ref{lem:ewwt2}, we get $\ell(xry) \geq \ell(ry)$. We have $\overline{\rho}_y(\rho_x(T_r)) = \overline{\rho}_y(T_{xr}) = q^{\ell(xr)-\ell(xry)}T_{xry}$. If $\ell(xry)  = \ell(ry) +1$, then  $\rho_x(\overline{\rho}_y(T_r)) =  \rho_x(q^{\ell(r)-\ell(ry)}T_{ry}) = q^{\ell(r)-\ell(ry)}T_{xry}$. If $\ell(xry)  = \ell(ry)$, then  $\rho_x(\overline{\rho}_y(T_r)) =  \rho_x(q^{\ell(r)-\ell(ry)}T_{ry}) = q^{\ell(r)-\ell(ry)+1}T_{xry}$. Thus, in both case, $\overline{\rho}_y(\rho_x(T_r)) = \rho_x(\overline{\rho}_y(T_r))$.\\
{\it Subcase~5:}~$x,y$ belong to~${S}$. If $\ell(xry) = \ell(r)-2$, then $\ell(xr) = \ell(ry) = \ell(r)-1$ and a calculation similar to~\cite[page~148 case (b)]{Hum} lied to~$\overline{\rho}_y(\rho_x(T_r)) =\rho_x( \overline{\rho}_y(T_r)) =  q^2T_{xry} +q(q-1)T_{xr} + q(q-1)T_{ry} + (q-1)^2T_{r}$. So, we consider the case~$\ell(xry) = \ell(r)-1$. If $\ell(xr) = \ell(r)$, then $xr = r$ and $xry = ry$. Therefore $\ell(ry)<\ell(r)$ and~$\overline{\rho}_y(\rho_x(T_r)) =\rho_x( \overline{\rho}_y(T_r)) = q(q-1)T_{xr}+q^2T_{xry}$. Now, consider the case $\ell(xr) = \ell(r)-1$. If $\ell(ry) = \ell(r)$, then~$\overline{\rho}_y(\rho_x(T_r)) =\rho_x( \overline{\rho}_y(T_r)) = q(q-1)T_{r}+q^2T_{xr}$; finally, if $\ell(ry) = \ell(r)-1$ then~$\overline{\rho}_y(\rho_x(T_r)) =\rho_x( \overline{\rho}_y(T_r)) = (q-1)^2T_{r}+q(q-1)T_{rt} + q^2T_{xry}$.
\end{proof}

Once we have Lemma~\ref{lem:commut}, we can almost repeat the argument of~\cite[Sec~7.3]{Hum} to prove Theorem~\ref{Th:genThe}. 

\begin{Lem} Let $\mathcal{L}$ be the sub-algebra of~$\mathcal{E}$ generated the $\rho_x$ for $x$ in $R$. The map~$\varphi$ from~$\mathcal{L}$ to $V$ which sends $\rho$ to $\rho(T_1)$ is an isomorphism of $A$-modules. \label{lem:iso}
\end{Lem}
\begin{proof} This is clear that~$\varphi$ is a morphism of $A$-modules. Let $r$ belong to~$R$, and let~$x_1\cdots x_k$ be a minimal word representative. Then by definition of the maps~$\rho_{x_i}$, we have $T_r = \varphi (\rho_{x_1}\cdots \rho_{x_k})$. Therefore, $\varphi$ is surjective. Assume $\varphi(\rho) = 0$ for some $\rho$ in~$\mathcal{L}$. Consider~$r$ and $x_1\cdots x_k$ as before, such that $k$ is minimal. We prove by induction on $k$ that $\rho(T_r) = 0$. For~$k = 0$, that is $r = 1$, this is true by hypothesis. The word~$x_1\cdots x_{k-1}$ is a minimal word representative of some element~$r'$. By induction hypothesis, we have $\rho(T_{r'}) = 0$. It follows $\rho(T_r) = \rho(T_{r'x_k}) = \rho(\overline{\rho}_{x_k}(T_{r'})) = \overline{\rho}_{x_m}(\rho(T_{r'})) = \overline{\rho}_{x_m}(0) = 0$.
\end{proof}

\begin{proof}[Proof of Theorem~\ref{Th:genThe}] Consider the notation of Lemma~\ref{lem:iso}. Assume $r$ belongs to $R$ and $x_1\cdots x_k$ is a minimal word representative of~$r$. Iterating the first defining relation in Theorem~\ref{Th:genThe}, we get $T_r = T_{x_1}\cdots T_{x_k}$. The unicity follows. Since $\varphi$ is an isomorphism, the endomorphism~$\rho_{r} = \rho_{x_1}\cdots\rho_{x_k}$ does not depend on the minimal word representing~$x_1\cdots x_k$, and the set~$\{\rho_r \mid r\in R\}$ is a free $A$-basis for~$\mathcal{L}$ with $\varphi(\rho_r) = \rho_r(T_1) = T_r$. Moreover, we can transfer the $A$-algebra structure of~$\mathcal{L}$ to $V$ using the isomorphism~$\varphi$. It remains to verify that the structure constants of the obtained $A$-algebra are the one stated in the theorem. Let $x$ belongs to~$S\cup \crlo$ and $r$ in~$R$. If $\ell(xr) = \ell(x)+\ell(r)$ and $\omega$ is a minimal word representative of $r$, then $x\omega$ is clearly a minimal word representative of~$xr$. Therefore~$\rho_x\rho_r(T_1) = \rho_x(T_r) = T_{xr} = \rho_{xr}(T_1)$. Therefore, $\rho_x\rho_r = \rho_{xr}$, and~$T_xT_r = T_{xr}$. Assume $x$ lies in $\crlo$ and $\ell(xr) < \ell(r)$. Then~$\rho_x\rho_r(T_1) = \rho_x(T_r) = q^{\ell(r)-\ell(xr)}T_{xr} = q^{\ell(r)-\ell(xr)} \rho_{xr}(T_1)$. We get $\rho_x\rho_r = q^{\ell(r)-\ell(xr)}\rho_{xr}$ and~$T_xT_r = q^{\ell(r)-\ell(xr)}T_{xr}$.  Assume $x$ lies in $S$. If $\ell(xr) = \ell(r)$, then~$\rho_x\rho_r(T_1) = \rho_x(T_r) = qT_{xr} = q\rho_{xr}(T_1)$ and $T_xT_r = qT_{rx}$. Finally, consider the case $\ell(xr) = \ell(r) -1$. One has ~$\rho_x\rho_r(T_1) = \rho_x(T_r) = (q-1)T_r + qT_{xr} = (q-1)\rho_r(T_1) + q\rho_{xr}(T_1) = ((q-1)\rho_r + q\rho_{xr})(T_1)$. Therefore,~$\rho_x\rho_r = (q-1)\rho_r + q\rho_{xr}$ and~$T_xT_r = (q-1)T_r + qT_{xr}$.
\end{proof}

\begin{definition} Let~$q$ be an indeterminate and set~$A = \mathbb{Z}[q]$. The \emph{generic Hecke algebra}~$\H(R)$ of the generalised Renner monoid~$R$ is the $A$-algebra described in Theorem~\ref{Th:genThe}. \label{def:genHechalg}
\end{definition}

\begin{Cor} The generic Hecke algebra~$\H(R)$ of~$R$ admits the following $\mathbb{Z}[q]$-algebra presentation: the generators are $T_x$ for $x$ in $S\cup\crlo$; the defining relations are\begin{center}   
\begin{tabular}{lll}
(HEC1)&$T^2_s = (q-1)T_1+qT_{s}$,&$s\in S$;\\
(HEC2)&$|T_s,T_t\rangle^m = |T_t,T_s\rangle^m$,&$(\{s,t\},m)\in \mathcal{E}(\Gamma)$;\\
(HEC3)&$T_sT_e = T_eT_s$,&  $e\in \crlo$, $s\in \lambda^\star(e)$;\\
(HEC4)&$T_sT_e = T_eT_s = qT_e$,& $e\in \crlo$, $s\in \lambda_\star(e)$;\\
(HEC5)&$T_eT_wT_f = q^{\ell(w)}T_{e\wedge_wf}$,& $e,f\in \crlo$, $w\in \textrm{Red}(e,f)$. 
\end{tabular} \end{center}\label{cor:preshHR}
\end{Cor}

In the special case of the rook monoid (see Example~\ref{exerook} below), we recover the presentation obtained in~\cite{God}.
\begin{proof}  Consider the presentation of~$\H(R)$ given in Theorem~\ref{Th:genThe}. Then Relations (HEC1)---(HEC5) clearly hold in~$\H(R)$. For instance~$|T_s,T_t\rangle^m  = T_{|s,t\rangle^m} = T_{|t,s\rangle^m} =  |T_t,T_s\rangle^m$.  Conversely, consider the algebra~$\H$ defined by the presentation given in the corollary. We claim that for two minimal word representatives~$\omega_1 = x_1\cdots x_k$ and $\omega_2 = y_1\cdots y_k$ on~$S\cup\crlo$ that represent the same element~$r$ in~$R$,  we have $T_{x_1}\cdots T_{x_k} = T_{y_1}\cdots T_{y_k}$. Indeed, it follows from Corollary~\ref{cor:pasaugleng} that we can transform  $T_{x_1}\cdots T_{x_k}$ into $T_{y_1}\cdots T_{y_k}$ by using $(HEC2)$, $(HEC3)$ and $(HEC5)$. So we set $T_r = T_{x_1}\cdots T_{x_k}$ in~$\H$. 
If $(w_1,e,w_2)$ is the normal decomposition of $r$ we have $T_r = T_{w_1}T_eT_{w_2}$. Now, we deduce that the defining relations of~$\H(R)$ given in Theorem~\ref{Th:genThe} hold in~$\H$ using lemma~\ref{lem:ewwt1} and~\ref{lem:ewwt2}. If $\ell(xr) = \ell(x)+\ell(r)$ and $x_1\cdots x_k$ is a minimal word representative of $r$, then $xx_1\cdots x_k$ is a minimal word representative of $xr$ and $T_{xr} = T_xT_{x_1}\cdots T_{x_k} = T_xT_r$. If $x$ belong to $S$ and $\ell(xr) = \ell(r) - 1$, then $T_xT_r =T_xT_{w_1}T_eT_{w_2} = ((q-1)T_{w_1} + qT_{xw_1}) T_eT_{w_2} = (q-1)T_r +qT_{xw_1}$. Here we use that Relations $(HEC1)$ and $(HEC2)$ implies $T_{w} = (q-1)T_{w} + qT_{xw}$ when $w$ belongs to $W$ such that $\ell(xw) = \ell(w) - 1$ ({\it cf.} \cite[Sec.~7]{Hum}). If $x$ belongs to $S$ and $\ell(xr) = \ell(r)$, then by Lemma~\ref{lem:ewwt1}, there exists~$u$ in~$\lambda_\star(e)$ such that $xw_1 = w_1u$, and~$\ell(xw_1) = \ell(w_1)+1$. It follows that~$T_xT_r = T_xT_{w_1}T_eT_{w_2} = T_{xw_1}T_eT_{w_2} = T_{w_1}T_uT_eT_{w_2} = qT_{w_1}T_eT_{w_2} = qT_r$. Finally, assume $x$ belongs to~$\crlo$ and $\ell(xr)<\ell(r)$. Write  $w_1 = w'''_1w''_1w'_1$ such that $\ell(w_1) = \ell(w'''_1)+\ell(w''_1)+\ell(w'_1)$ with~$w'''_1$ in~$W_\star(x)$,~$w''_1$ in~$W^\star(x)$ and~$w'_1$ in $\textrm{Red}(x,e)$. We have~$T_xT_r = T_xT_{w_1}T_eT_{w_2} = T_xT_{w'''_1}T_{w''_1}T_{w'_1}T_eT_{w_2} = q^{\ell(w'''_1)}T_{x} T_{w''_1}T_{w'_1}T_eT_{w_2} = q^{\ell(w'''_1)} T_{w''_1}T_{x}T_{w'_1}T_eT_{w_2}$.  We get~$T_xT_r  = q^{\ell(w'''_1)+\ell(w'_1)} T_{w''_1}T_{x\wedge_{w'_1}e}T_{w_2}$. We can decompose~$w''_1$ and~$w_2$ such that $w''_1 = v'_1v''_1$ and $w_2 = v''_2v'_2$ where $v''_1,v''_2$ belong to~$W_\star(x\wedge_{w'_1}e)$, $v'_1$ belongs to $\DrTe{x\wedge_{w'_1}e}$ and $v'_2$ belongs to $\GaTe{x\wedge_{w'_1}e}$. We have~$\ell(xr) = \ell(v'_1)+\ell(v'_2)$ and $v'_1(x\wedge_{w'_1}e)v'_2$ is a minimal word representative of $xr$. Hence, $T_xT_r = q^{\ell(w'''_1)+\ell(w'_1)+\ell(v''_1)+\ell(v''_2)} T_{v'_1}T_{x\wedge_{w'_1}e}T_{v'_2} = q^{\ell(x)-\ell(xr)} T_{xr}$.    
\end{proof}

\begin{Rem}\label{Rem:rem2}
(i) For $e,f$ in~$\crlo$,  we set $$\textrm{Red}_\star(e,f) = \textrm{Red}(e,f)\bigcap W_{\cap_{h> e}\lambda(h)} \bigcap W_{\cap_{h> f}\lambda(h)}.$$ It is not difficult to see that in Relations~$(HEC5)$ of the presentation stated in Corollary~\ref{cor:preshHR}, we can assume $w$ belongs to $\textrm{Red}_\star(e,f)$ ({\it cf.} the proof of ~\cite[Theorem~0.1]{God2}).\\(ii) In $\H(R)$ the following relations hold : 
\begin{center}\begin{tabular}{ll}$T_rT_x = T_{xr}$,&if $x\in S$ and $\ell(rx) =  \ell(r)+1$;\\
$T_rT_x = qT_{r}$,&if $x\in S$ and $\ell(rx) = \ell(r)$;\\$T_rT_x = (q-1)T_r+qT_{rx}$,&if $x\in S$ and $\ell(rx) = \ell(r) - 1;$\\ $T_rT_x = q^{\ell(r)-\ell(rx)}T_{rx}$,&if $x\in \crlo$.
\\\end{tabular}\end{center} 
This can be deduced directly from Theorem~\ref{Th:genThe}, but this is an immediate consequence of Corollary~\ref{cor:preshHR} since the defining relations~$(HEC1)-(HEC5)$ have a right-left symmetry.
\end{Rem}

\section{Iwahori-Hecke algebra of finite reductive monoids}
Here, we first recall basic results on Algebraic Monoid Theory, then we introduce the notion of an Iwahori-Hecke algebra in the general framework of Monoid Theory, we recall some basic properties and explain why this Iwahori-Hecke algebra is interesting. Finally, we turn to finite reductive monoids and prove that the Iwahori-Hecke algebra of such monoids is related to the generic Hecke algebra of the associated Renner monoid. As a consequence, we prove Theorems~\ref{THintro} and \ref{THintro2}. 
\label{sectionprincipale}
\subsection{Regular monoids and reductive groups}
\label{sousect}
We introduce here the basic definitions and 
notation on Algebraic Monoid Theory that we shall need in the sequel.  We fix an algebraically closed field~$\mathbb{K}$. We let~$M_n$ denote the set of all~$n\times n$ matrices over~$\mathbb{K}$, and by~$GL_n$ the set of all invertible matrices in~$M_n$. We refer to~\cite{Put,Ren,Sol1} for the general theory and proofs involving linear algebraic monoids and Renner monoids; we refer to~\cite{Hum} for an introduction to Linear Algebraic Groups Theory. If $X$ is a subset of~$M_n$, we let $\overline{X}$ denote its closure for the Zariski topology. Recall that a semigroup~$M$ is said to have a \emph{zero element} if it contains an element~$0$ such that~$0\times x = x\times 0 = 0$ for every~$x$ in~$M$. 
 
\begin{definition}[Algebraic monoid] An \emph{algebraic monoid} is a submonoid of~$M_n$, for some positive integer~$n$, that is closed for the Zariski topology. An algebraic monoid is \emph{irreducible} if it is irreducible as a variety.\end{definition}

It is very easy to construct algebraic monoids. Indeed, the Zariski closure~$M = \overline{G}$ of any submonoid~$G$ of~$M_n$ is an algebraic monoid. The main example occurs when for $G$ one considers an algebraic subgroup of~$GL_n$. It turns out that in this case, the group~$G$ is the unit group of~$M$. Conversely, if $M$ is an algebraic monoid, then its unit group~$G(M)$ is an algebraic group. The monoid~$M_n$ is the seminal example of an algebraic monoid, and its unit group~$GL_n$ is the seminal example of an algebraic group.

The next result, which is the starting point of the theory, was obtained independently by Putcha and Renner in 1982. 
\begin{The} Let $M$ be an irreducible algebraic monoid with a zero element. Then $M$ is regular if and only if its unit group~$G(M)$ is reductive.  
\end{The}

\begin{definition}[Reductive monoid] A \emph{reductive monoid} is an irreducible algebraic monoid whose unit group is a reductive group. \end{definition}

\begin{definition}[Renner monoid] Let $M$ be a reductive monoid. The normaliser of a maximal torus~$T$ of $G(M)$ is denoted by $N_{G(M)}(T)$. The \emph{Renner monoid}~$\renn$ of~$M$ is the monoid~$\overline{N_{G(M)}(T)}/T$. 
\end{definition}
It is clear that~$\renn$ does not depend on the choice of the maximal torus of the algebraic group~$G(M)$.

\begin{Prop} Let $M$ be reductive monoid. Fix a maximal torus~$T$ of $G(M)$ and a Borel subgroup~$B$ of $G(M)$ that contains~$T$. The unit group of~$R(M)$ is the Weyl group~$W$ of $G(M)$. If $S$ is the standard generating set of~$W$ associated with the Borel $B$ and $\crl(B) = \{e\in E(\overline{T})\mid \forall b\in B,\ be = ebe \}$, then $(R(M),\crl(B),S)$ is a generalised Renner-Coxeter system such that~$R(M)$ is a generalised Renner monoid. Moreover, there is a canonical order preserving isomorphism of monoids between $E(R(M))$ and $E(\overline{T})$.\label{propRgenparEetW}
\end{Prop}

\begin{Exe} \label{exerook}Consider $M = M_n$. Choose the Borel subgroup~$\mathbb{B}$ of invertible upper triangular matrices and the maximal torus~$\mathbb{T}$ of invertible diagonal matrices. The Renner monoid is isomorphic to the monoid of matrices with at most one nonzero entry, that is equal to~$1$, in each row and each column. This monoid is called the rook monoid~$R_n$~\cite{Sol2}. Its unit group is the group of monomial matrices, which is isomorphic to the symmetric group~$S_n$.  Denote by~$e_i$  the diagonal matrix~$\left(\begin{array}{cccccc}Id_i&0\\0&0\end{array}\right)$ of rank~$i$. Then the set~$\crl(\mathbb{B})$  is~$\{e_0,\ldots, e_n\}$. One has $e_i\leq e_{i+1}$ for every index~$i$. One has~$\lambda_\star(e_i) = \{s_j\mid j > i\}$ and $\lambda^\star(e_i) = \{s_j\mid j < i\}$.\begin{figure}[ht]
\begin{picture}(250,75)
\put(18,0){\includegraphics[scale = 0.6]{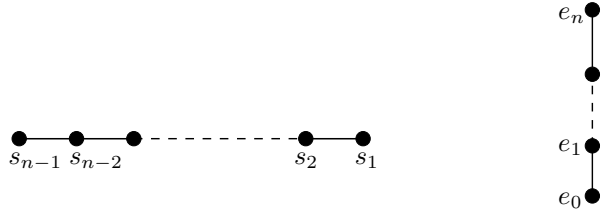}}
\put(148,15){$s_1$} \put(125,15){$s_2$}\put(40,15){$s_{n-2}$}\put(17,15){$s_{n-1}$}
\put(225,0){$e_0$} \put(225,20){$e_1$} \put(225,70){$e_n$}
\end{picture}
\caption{Coxeter graph~$\Gamma(S)$ and Hasse diagram~$\Lambda(\mathbb{B})$ for $M_n$.}\label{fig:hassecrlAn}
\end{figure}\\
Other examples can be found in~\cite{God2}. \end{Exe}

In the framework of algebraic monoids, Renner monoid plays the role of Weyl groups in Algebraic Group Theory. In particular we still have a \emph{Bruhat decomposition}: the monoid~$M$ is equal to the disjoint union~$\cup_{r\in R}BrB$. Moreover,  the product of double classes~$BrB$ is related to the length function that we introduce in Section~\ref{seclong}:

\begin{Prop} Let $M$ be a reductive monoid. Fix a maximal torus~$T$ of $G(M)$ and a Borel subgroup~$B$ of $G(M)$ that contains~$T$. Consider the generalised Renner-Coxeter system~$(\renn,\Lambda,S)$ of~$\renn$ defined in Proposition~\ref{propRgenparEetW}.\\(i) Let~$r$ lie in~$\renn$ and~$s$ lie in~$S$, then $$B s B r B = \left\{\begin{array}{ll}BrB,&\textrm{if } \ell(sr) = \ell(r);\\BsrB,&\textrm{if } \ell(sr) = \ell(r)+1;\\BsrB\cup BrB,&\textrm{if } \ell(sr) = \ell(r)-1.\end{array}\right.$$ \label{Pr:lienlongBB}
(ii) Let~$r$ lie in~$\renn$ and~$s$ lie in~$S$, then $$B r B s B = \left\{\begin{array}{ll}BrB,&\textrm{if } \ell(rs) = \ell(r);\\BrsB,&\textrm{if } \ell(rs) = \ell(r)+1;\\BrsB\cup BrB,&\textrm{if } \ell(rs) = \ell(r)-1.\end{array}\right.$$\\(iii) Let~$r$ lie in~$\renn$ and~$e$ lie in~$\Lambda$, then $$B e B r B = BerB \textrm{ and } BrBeB = BreB$$ \end{Prop}
\begin{proof} (i) is proved in~\cite[Prop.~0.2]{God2} in the case of irreducible regular monoid~$M$ with a zero element. Same arguments can be applied for any reductive monoids; let us deduced~(ii): by the remark following~\cite[Prop.~8.6]{Ren} we know that~$$B r B s B\subseteq BrB\cup BrsB$$  and, clearly, $B r B s B$ is a union of double classes. Hence, $B r B s B$  has to be equal to $BrB$, $BrsB$ are $BrB\cup BrsB$. If~$\ell(rs) = \ell(r)$ then $rs = r$ and we are done. if $\ell(rs) = \ell(r)+1$ and $r = x_1\cdots x_k$ is a minimal word representative of $r$ then $BrBsB = Bx_1B\cdots Bx_{k-1}Bx_kBsB = Bx_1B\cdots Bx_{k-1}Bx_ksB = \cdots = BrsB$. Finally, if $\ell(rs) = \ell(rs)-1$, and $x_1\cdots x_{k-1}s$ is a minimal word representative of $r$, then $BrBsB = Bx_1B\cdots Bx_{k-1}BsBsB = Bx_1B\cdots Bx_{k-1}B (B\cup BsB) = BrsB \cup BrB$. Let us proof (iii). Since $e$ belongs to $\Lambda$, $Be\subseteq eB$ \cite{Ren}. Thus, $BrBeB\subseteq BreB$. The inclusion~$BreB\subseteq BrBeB$ is trivial. Let us prove that $BeBrB = BerB$.
If $r = s_{i_1}\cdots s_{i_{\ell(r)}}$ belongs to the Weyl group~$W$, the results follows from~(ii) since for $\ell(es_{i_1}\cdots s_{i_j})\geq \ell(es_{i_1}\cdots s_{i_{j-1}})$. Therefore, we may assume that $r = w_1fw_2$ where $f$ lies in~$\crlo$ and $(w_1,f,w_2)$ is the normal decomposition of~$r$.
We can write $w_1 = v_1v_2v_3v_4$ with $v_1\in W_\star(e)$, $v_2\in W^\star(e)$,  $v_3\in \textrm{Red}(e,f)$, $v_4$ in $W^\star(f)$ and $\ell(w_1) = \ell(v_1)+\ell(v_2)+\ell(v_3)+\ell(v_4)$. Then $$\displaylines{BeBrB = BeBw_1fw_2B = BeBv_1v_2v_3fv_4w_2B = BeBv_1Bv_2v_3BfBv_4w_2B=\hfill\cr\hfill Bev_2v_3BfBv_4w_2B = Bv_2ev_3fBv_4w_2B = Bv_2(e\wedge_{v_3}f)Bv_4w_2B.}$$ Write $v_4w_2 = v_5v_6v_7$ such that $\ell(v_4w_2) = \ell(v_5)+\ell(v_6)+\ell(v_7)$ and $v_5\in W_\star(e\wedge_{v_3}f)$, $v_6\in W^\star(e\wedge_{v_3}f)$, $v_7\in \Gae{e\wedge_{v_3}f}$. Then  $BeBrB = Bv_2(e\wedge_{v_3}f)Bv_6v_7B$. We claim that~$\ell(er) = \ell(v_2(e\wedge_{v_3}f)v_6v_7) = \ell(v_2(e\wedge_{v_3}f)) + \ell(v_6v_7)$, which implies $BeBrB = Bv_2(e\wedge_{v_3}f)v_6v_7B = BerB$ by (ii). If it was not the case, By Lemma~\ref{lem:ewwt1}~(iii), $v_6v_7 = uv_8$ with $u\in \lambda^\star(e\wedge_{v_3}f)$, $\ell(v_6v_7) = \ell(v_8)+1$ and $\ell(v_2u) = \ell(v_2)-1$. But $\lambda^\star(e\wedge_{v_3}f)\subseteq \lambda^\star(f)$, $uv_5 = v_5u$ and $uv_2 = v_2u$ since $v_2$ lies in $W_\star(e\wedge_{v_3}f)$. Therefore, this leads to $r = w_1ew_2 = v_1v_2v_3v_4fw_2 = v_1v_2v_3fv_4w_2 = v_1v_2v_3fv_5uv_8 = v_1v_2uv_3fv_5v_8$. But this is impossible since $$\displaylines{\ell(r) = \ell(v_1v_2uv_3fv_5v_8)\leq \ell(v_1)+\ell(v_2u)+\ell(v_3)+\ell(v_5)+\ell(v_8)=\hfill\cr\hfill \ell(v_1)+\ell(v_2)-1+\ell(v_3)+\ell(v_5)+\ell(v_8)< \ell(w_1)+\ell(w_2) = \ell(r).}$$   
\end{proof}

\subsection{Iwahori-Hecke algebra}
We introduce here the notion of a Iwahori-Hecke algebra in the general framework of Monoid Theory. The equivalent notion in the context of Group Theory is well-known~(\cite[Sec.~8.4]{Gec} for instance). There is no difficulty to translate the notion from Group Theory to Monoid Theory. The point is to verify that definitions and proofs can be written without using the existence of inverse elements. This is not the case for the whole theory (see  Remarks~\ref{rem:1} and~\ref{rem:2} below) but the main results still hold as far as one considers the Iwahori-Hecke algebra associated with a subgroup. We have no find general references for Iwahori-Hecke Algebra of a monoid. This is why we start with an introduction to these notions with included proof. 
\label{sect:monhecke}
 
For all this section, we assume $M$ is a finite monoid. We let $G$ denote its unit group and we fix a subgroup~$H$ of $G$. 
We let~$\C[M]$ denote the monoid algebra of~$M$. An element of~$\C[M]$ has the form $\sum_{x\in M}\lambda_x x$ where the $\lambda_x$ belong to~$\C$. We set $$\ee = \frac{1}{|H|}{\sum_{h\in H}}h$$ in~$\C[M]$. All the considered algebras are unit associative algebras, and all modules are left modules. We begin with two easy lemma whose proofs are left to the reader. \\
\begin{Lem}\label{lem:basic1}
Consider the $\C$-algebra~$\C^M$ of linear maps from $M$ to~$\C$ where the product is the convolution product~$\star$, defined by $$f\star g (x) = \sum_{y,z\in M, yz = x}f(y)g(z).$$ There is a canonical isomorphism of $\C$-algebra from~$\C[M]$ to $\C^M$ which sends~$X = \sum_{x\in M}\lambda_x x$ to the map~$\overline{X}: x\mapsto\lambda_x$.  
\end{Lem}
The following lemma is immediate. We left the proof to the reader.
\begin{Lem}\label{lem:basic2}
(i) $\ee^2 = \ee$, and for every~$h$ in $H$ one has $h\ee = \ee h = h$.\\(ii) $\C[M]\ee$ and $\C[M\slash H]$ are isomorphic as~$\C[M]$-modules and as~$\C$-vector spaces. 
\end{Lem}
\begin{Rem}We remark that Lemma~\ref{lem:basic2} is no more true in general if we only assume~$H$ is a submonoid of $M$. Indeed, $\ee$ is not necessarily an idempotent.\label{rem:1}\end{Rem} 
\begin{Prop} \label{Prop:algisom}There is a canonical isomorphism between the following $\C$-algebras:\\(a) the subalgebra of $\C^M$ whose elements are the linear maps which are constant on the double-classes $H\backslash M\slash H$;\\
(b) the algebra~$\ee \C[M] \ee$;\\
(c) the algebra~$\left(End_{\C[M]}(\C[M\slash H])\right)^{op}$ of endomorphisms of $\C[M\slash H]$  considered as a $\C[M]$-module (for the opposite product).
 \end{Prop}
\begin{proof}
The second and third algebras are isomorphic by \cite[Lemma~3.19]{CurRei}.  This is clear that $\ee X \ee = X$ if and only if $X$ belongs to~$\ee \C[M] \ee$. Consider the notation of Lemma~\ref{lem:basic1}. Denote by $Hx_1,\ldots, Hx_k$  the left classes of $M$ modulo the subgroup~$H$. Let $X = \sum_{x\in M}\lambda_x x$ belong to~$\C[M]$. Then $$\ee X = \frac{1}{|H|}\sum_{i = 1}^k\sum_{x\in H\!x_i}\sum_{h\in H} \lambda_{x} hx = \sum_{i = 1}^k\sum_{x\in H\!x_i}\left(\frac{1}{|H|}\sum_{y\in Hx_i} \alpha_{y,x}\lambda_{y} \right) x$$ where $\alpha_{y,x} = \#\{h\in H\mid hy = x\}$. If $M$ is a group, then $\alpha(y,x) = 1$ for every $y,x$ in $Hx_i$. In the general case one has $\alpha(y,x) = \frac{|H|}{|Hx_i|}$ because $H$ is a group. Therefore, $\ee X = \sum_{x\in M}\left(\frac{1}{|Hx|}\sum_{y\in Hx}\lambda_{y} \right) x$, and  $\ee X = X$ if and only if $\overline{X}$ is constant on each left class. by a similar computation, $X \ee = X$ if and only if $\overline{X}$ is constant on each right class. Therefore $\ee X \ee = X$ if and only if $\overline{X}$ is constant on each double class. 
\end{proof}
\begin{Rem} The isomorphism between~$\left(End_{\C[M]}(\C[M\slash H])\right)^{op}$ and~$\ee \C[M] \ee$ is given by $f\mapsto \ee f(\ee)\ee$ for every endomorphism~$f$.
\end{Rem}

Following Solomon~\cite{Sol} and Putcha~\cite{Put}, who consider the case of finite reductive monoids, we introduce the Iwahori-Hecke algebra~$\H(M,H)$:
\begin{definition}[Iwahori-Hecke algebra] Let $M$ be a finite monoid, and assume $H$ is a subgroup of $M$. Let~$\ee = \frac{1}{|H|}\sum_{h\in H}h$ in~$\C[M]$. We define the~\emph{Iwahori-Hecke algebra~$\H(M,H)$ of $M$ relatively to $H$} to be the algebra~$\ee \C[M] \ee$.  
\end{definition}

It is immediate that for every $\C[M]$-module~$N$, we get an induced structure of left $\H(M,H)$-module on~$\ee N$. Proposition~\ref{Prop:algisom} explains why the Hecke algebra is interesting. Another motivation for such a definition is the following result. 
\begin{Prop} Assume~$\C[M]$ is semisimple.\\(i) The Hecke algebra~$\H(M,H)$ is semisimple.\\ (ii) The map $N\mapsto \ee N$ induced a one-to-one correspondence between the set of simple~$\C[M]$-modules in the induced~$\C[M]$-module $\C[M]\ee = \C[M]\otimes_{\C[H]}\C[H]$ and the set of isomorphic classes of simple $\H(M,H)$-modules. Furthermore, the multiplicity of $N$ in $\C[M]\ee$ is equal to the dimension of the $\H(M,H)$ module~$\ee N$ considered as a $\C$-vector space. 
 \end{Prop}

Note that this is known by~\cite{OkPu} that~$\C[M]$ is semisimple for abstract finite monoids of Lie type ({\it cf.} Example~\ref{exe:ex2}), and therefore for finite reductive monoids. 
\begin{proof}
Since~$\C[M]$ is semisimple, the algebra~$\ee \C[M] \ee$ is semisimple. Assume $N$ is a simple~$\C[M]$ module and let~$f$ belong to $Hom_{\C[M]}(\C[M]\ee,N)$. For every $x$ in $\C[M]\ee$ one has $f(x) = f(x\ee) = xf(\ee)$. If we consider $x = \ee$, we get that $f(\ee)$ belongs to~$\ee N$. Moreover, it follows that the map~$f\mapsto f(\ee)$ from $Hom_{\C[M]}(\C[M]\ee,N)$ to $\ee N$ is $\C$-linear and one-to-one. Thus $dim_\C(\ee N)$ is equal to $dim(Hom_{\C[M]}(\C[M]\ee,N))$, that is to the multiplicity of $N$ in~$\C[M]\ee$. Now write $\C[M]\ee = \oplus_i M_i$ where the $M_i$ are simple $\C[M]$-modules. Then $\ee\C[M]\ee = \oplus_i \ee M_i$ and each~$\ee M_i$ is a non-trivial simple $\H(M,H)$-modules: its $\C$-dimension is at least one, and for $m$ in $M_i$ such that $\ee m\neq 0$ one has $\H(M,H)\ee m = \ee\C[M]\ee m = \ee M_i$ since $M_i$ is a simple $\C[M]$-module.       
\end{proof}
By Proposition~\ref{Prop:algisom}, this is immediate to obtain a $\C$-basis of $\H(M,H)$:.
\begin{Prop} Let $\{D_1,\cdots, D_\ell\}$ be the set of double classes of $M$ modulo~$H$. We fix some arbitrary non-zero complex numbers~$a_1,\cdots a_\ell$, and  we set $\displaystyle X_i = a_i\sum_{x\in D_i}x$ for $i$ in $\{1,\cdots,\ell\}$. Then the $X_i$ form a $\C$-basis for $\H(M,H)$. If we write $\displaystyle X_iX_j = \sum_{k = 1}^\ell \mu(i,j,k)X_k$, then $\mu(i,j,k) = \frac{a_ia_j}{a_k}\#\{(x,y)\in D_i\times D_j\mid xy = x_k\}$ where $x_k$ is an arbitrary fixed element of $D_k$. \end{Prop}
\begin{proof} The first part is clear. The second part come from the fact that $H$ is a group: we can write~$X_iX_j = \sum_{k = 1}^\ell \sum_{z\in D_k} \alpha(i,j,z) z$ where $\alpha(i,j,z)  = \#\{(x,y)\in D_i\times D_j\mid xy = z\}$. But if $z$ belongs to $D_k$, then $\alpha(i,j,z) = \alpha(i,j,x_k)$. Indeed, if $z = h_1x_kh_2$ then the map $(x,y)\mapsto (h_1x,yh_2)$ is one-to-one from $\{(x,y)\in D_i\times D_j\mid xy = x_k\}$ onto $\{(x,y)\in D_i\times D_j\mid xy = z\}$.      
\end{proof}

As explained in~\cite[Sec~4]{Sol} and in~\cite[Sec~2]{Put}, an important issue is to determined the structure constants~$\mu_{i,j,k}$ and, if possible, to suitably choose the $a_i$ so that the $\mathbb{Z}$-module generated by the~$a_i X_i$ becomes a $\mathbb{Z}$-subalgebra of $\H(M,B)$, in other words, so that the structure constants~$\mu_{i,j,k}$ belong to~$\mathbb{Z}$.     
\begin{Rem}\label{rem:2} Let $\varphi$ belong to $End_{\C}(\C[M\slash H])$. Define $\dot{\varphi} : M\slash H\times M\slash H\to \mathbb{C}$ by $\varphi(xH) = \sum_{yH\in M\slash H}\dot{\varphi}(yH,xH)yH$. If $M$ is a group, it turns out that $\varphi$ belongs to~$End_{\C[M]}(\C[M\slash H])$, that is to $\H(M,H)$, if and only if $\dot{\varphi}$ is constant on the orbits of~$M$ on $M\slash H\times M\slash H$~\cite[Sec~8.4]{Gec}, which are naturally related to the double classes $HxH$ when $M$ is a group. This is no more true if we only assume $M$ is a monoid. One can verify that in the general case, $\varphi$ belongs to~$End_{\C[M]}(\C[M\slash H])$ if and only for every $xH$ and $yH$ in $M/H$ and every $g$ in $M$, one has $\dot{\varphi}(yH,gxH) = 0$ if $yH\cap gM$ is empty, and $$\dot{\varphi}(gyH,gxH) = \frac{1}{|\mathcal{C}_g(yH)|}\sum_{zH\in \mathcal{C}_g(yH)}\dot{\varphi}(zH,xH)$$ where $\mathcal{C}_g(yH) = \{zH\mid gzH = gyH\}$. If $M$ is a group then $yH\cap gM$ is never empty, and  $\mathcal{C}_g(yH) = \{yH\}$. 
\end{Rem}
\subsection{Finite reductive monoids}
\label{sousectionrappelamt}
We can now turn to the proof of Theorems~\ref{THintro} and~\ref{THintro2}. Let us recall the definition of finite reductive monoids~\cite{Ren2}, which is in the spirit of the definition of finite reductive groups~\cite{Ste}.
\begin{definition}[finite reductive monoid] Let~$\underline{M}$ be a reductive monoid defined over $\overline{\mathbb{F}}_q$. A finite submonoid~$M$ of~$\underline{M}$ is a \emph{finite reductive monoid} if there exists a surjective endomorphism of algebraic monoid~$\sigma : \underline{M}\to \underline{M}$ such that $$M = \{x\in \underline{M}\mid \sigma(x) = x\}.$$\label{def:finredmon} \end{definition}

\begin{Exe}
Consider a reductive monoid~$\underline{M}$ over $\overline{\mathbb{F}}_q$. The  finite reductive monoid~$M$ associated with the map~$(x_{i,j})\mapsto (x_{i,j}^q)$ is~$M_n(\mathbb{F}_q)$. See~\cite{Sol} for more details. 
\end{Exe}

Finite reductive monoids are special cases of abstract finite monoids of Lie type~\cite{Put3}, and their unit groups are finite groups of Lie type. Therefore, they are \emph{groups with a BN pair} and possess Borel subgroups and a generalised Renner monoid~$R$ ({\it cf.} Example~\ref{exe:ex2}). As a consequence, we can associate with $M$ a generic Hecke algebra~$\H(R)$ as defined in Section~2, and a Iwahori-Hecke algebra as defined in Section~\ref{sect:monhecke}. Our objective is to prove Theorem~\ref{THintro2}, which explains how these two notions are related.

\begin{Not} Assume~$M$ is a finite reductive monoid over $\overline{\mathbb{F}}_q$, and consider the notation of Definition~\ref{def:finredmon}. There exists a maximal torus~$\underline{T}$ of~$G(\underline{M})$ and a Borel subgroup~$\underline{B}$ of $\underline{G} = G(\underline{M})$ that contains $\underline{T}$ such that $\sigma(\underline{T}) = \underline{T}$ and $\sigma(\underline{B}) = \underline{B}$ \cite{Ste,Ren2}. Moreover, $\sigma(N_{\underline{G}}(\underline{T})) = N_{\underline{G}}(\underline{T})$. Let $\underline{R}$ be the Renner monoid associated with~$\underline{M}$, and $\underline{W}$ be its unit group. Then $\sigma$ induces an isomorphism~$\sigma : \underline{R} \to \underline{R}$. We set $$G = \{b\in \underline{G}\mid \sigma(g) = g\}$$ $$B = \{b\in \underline{B}\mid \sigma(b) = b\}$$ $$T = \{t\in \underline{T}\mid \sigma(t) = t\}$$ $$W = \{w\in \underline{W}\mid \sigma(w) = w\}$$ $$R = \{r\in \underline{R}\mid \sigma(r) = r\}$$ $$\crl = \{e\in \underline{\crl}\mid \sigma(e) = e\}$$\label{not:fintredmon}
\end{Not}
\begin{Prop}\cite{Ren2,Ste} Consider Notation~\ref{not:fintredmon}. The group~$G$ is the unit group of~$M$, and $B$ is a Borel subgroup of~$G$ with maximal torus~$T$. The Renner monoid of $M$ is~$R$. The unit group of $R$ is $W$, and $\Lambda$ is the cross section lattice of~$R$ associated with~$B$. Denote by $\underline{S}$  the canonical generating set of~$\underline{W}$ associated with~$\underline{T}$ and~$\underline{B}$. For a conjugated class~$X$ of elements of~$S$ under~$\sigma$, we let~$\Delta_X$ denote the greatest element of~$W_X$. Let $S$ be the set of all $\Delta_X$. Then~$(W,S)$ is a Coxeter system, and $(R,\Lambda,S)$ is a generalised Renner-Coxeter system. Moreover, we have a disjoint union \emph{Bruhat decomposition}~$M = \cup_{r\in R}BrB$. \label{Prop:redtofinred} 
\end{Prop}

From the Bruhat decomposition of~$M$, we deduce for every $r$ in $R$ that $$BrB = \{x\in \underline{B}r\underline{B}\mid \sigma(x) = x\}.$$ It is immediate that for $e$ in~$\Lambda$ one has $\sigma(\underline{\lambda}(e)) = \underline{\lambda}(e)$  and $\sigma(\underline{\lambda}_\star(e)) = \underline{\lambda}_\star(e)$ in $\underline{R}$, with obvious notation. Therefore, $\omega_X$ belongs to $\lambda(e)$ in $R$  ({\it resp.} to~$\lambda_\star(e)$) if and only if $X$ is included in~$\underline{\lambda}_(e)$ ({\it resp.} to~$\underline{\lambda}_\star(e)$) in $\underline{R}$. 
\begin{Lem} Consider Notation~\ref{not:fintredmon}. Denote by $\ell$ the length function on~$R$.\\ (i) Let~$r$ lie in~$R$ and~$s$ lie in~$S$. Then $$B s B r B = \left\{\begin{array}{ll}BrB,&\textrm{if } \ell(sr) = \ell(r);\\BsrB,&\textrm{if } \ell(sr) = \ell(r)+1;\\BsrB\cup BrB,&\textrm{if } \ell(sr) = \ell(r)-1.\end{array}\right.$$
 (ii) Let~$r$ lie in~$R$ and~$s$ lie in~$S$. Then $$B r B s B = \left\{\begin{array}{ll}BrB,&\textrm{if } \ell(rs) = \ell(r);\\BrsB,&\textrm{if } \ell(rs) = \ell(r)+1;\\BrsB\cup BrB,&\textrm{if } \ell(rs) = \ell(r)-1.\end{array}\right.$$
(iii) Let $e$ lie in $\crlo$ and $r$ lie in $R$. Then $$B e B r B = BerB.$$\end{Lem}

\begin{proof} The result follows from Proposition~\ref{Pr:lienlongBB}. (i) Denote by $\underline{\ell}$ the length function on~$\underline{R}$. Let $r$ lie in $R$ and $\Delta_X$ lie in $S$ ({\it cf.} Proposition~\ref{Prop:redtofinred}). Fix a minimal representative word~$x_1\cdots x_k$  on~$\underline{S}$ of $\Delta_X$. Using the map~$\sigma$, we deduce that there is three possibilities:\\ (a) $\forall x \in X$, $\underline{\ell}(xr) = \underline{\ell}(r)+1$.\\ In this case, $\underline{\ell}(\omega_Xr) = \underline{\ell}(r)+\underline{\ell}(\omega_X)$, $\underline{B} \omega_X \underline{B} r \underline{B} = \underline{B} \omega_Xr \underline{B}$ and $\ell(\omega_Xr) = \ell(r)+1$. Therefore, $B \omega_X B r B \subseteq \{x\in\underline{B} \omega_Xr \underline{B}\mid \sigma(x) = x\} = B \omega_Xr B$. But $B \omega_X B r B$ is an union of double classes~$ByB$. Then the latter inclusion has to be an equality.\\  (b) $\forall x \in X$, $\underline{\ell}(xr) = \underline{\ell}(r)$.\\
In this case $\omega_Xr = r $, and in particular $\underline{\ell}(\omega_Xr) = \underline{\ell}(r)$, $\underline{B} \omega_X \underline{B} r \underline{B} = \underline{B} r \underline{B}$ and $\ell(\omega_Xr) = \ell(r)$. It follows that~$B \omega_X B r B = B r B$ as in the previous case.\\
  (c) $\forall x \in X$, $\underline{\ell}(xr) = \underline{\ell}(r)-1$.\\ In this case, $\underline{\ell}(\omega_Xr) = \underline{\ell}(r)-\underline{\ell}(\omega_X)$, $\ell(\omega_Xr) = \ell(r)-1$  and $\underline{B} \omega_X \underline{B} r \underline{B} = \bigcup_{v}\underline{B} vr \underline{B}$, where~$v$ ranges over all the elements $x_{i_1}\cdots x_{i_j}$ with $1\leq i_1<\cdots< i_j\leq k$ and $0\leq j\leq k$. But for such an element~$v$ of~$\underline{R}$, the set $\{x\in \underline{B} vr \underline{B}\mid \sigma(x) = x\}$ is empty, except if $vr$ belongs to~$R$, that is $v = 1$ or $v = \omega_X$.  Therefore, $\{x\in \underline{B} \omega_X \underline{B} r \underline{B}\mid \sigma(x) = x\} = B \omega_Xr B\cup B r B$. But~$\underline{B} \omega_X \underline{B} r \underline{B} = \bigcup_{b\in \underline{B}}\underline{B}\omega_Xbr \underline{B}$. We deduce that~$$M\cap \underline{B} \omega_X \underline{B} r \underline{B} = \bigcup_{b\in \underline{B}}M\cap \underline{B}\omega_Xbr \underline{B} = \bigcup_{b\in B}B\omega_Xbr B = B\omega_XBr B.$$(ii) the proof is similar to~(ii). \\(iii) $BeBrB$ is included in $\{x\in \underline{B} er \underline{B}\mid \sigma(x) = x\} = BerB$. But $BeBrB$ is an union of double classes $ByB$. Therefore, $BeBrB = BerB$. 
\end{proof}


We are now ready to prove Theorem~\ref{THintro2}.
\begin{proof}[Proof of Theorem~\ref{THintro2}] By Theorem~\ref{Th:genThe} and  Definition~\ref{def:genHechalg}, $\mathbb{C}\otimes_{\mathbb{Z}}\H_q(R)$ is the unique $\mathbb{C}$-algebra such that the relations stated in Theorem~\ref{Th:genThe} hold. But, by Section~\ref{sect:monhecke}, $\H(M,B)$ is a $\mathbb{C}$-algebra over the free $\mathbb{C}$-module with basis $\sum_{x\in BrB}x$, for $r\in R$. We set $$T_r = \frac{q^{\ell(r)}}{|BrB|} \sum_{x\in BrB}x$$ in~$\H(M,B)$. We are going to prove that the relations stated in Theorem~\ref{Th:genThe} hold
in~$\H(M,B)$ for the basis $T_r$, $r\in R$. The main arguments are like in~\cite[Sec.~4]{Sol}. Denote by~$$\pi: \H(M,B)\to \mathbb{C}$$ the restriction of the one-dimensional representation from~$\mathbb{C}[M]\to\mathbb{C}$ that sends every $g$ in $M$ to $1$. We have $\pi (T_r) = \pi(\frac{q^{\ell(r)}}{|BrB|} \sum_{x\in BrB}x) = q^{\ell(r)}$.  Let $r_1,r_2,r_3$ lie in $R$ such that $Br_1B Br_2B = Br_3B$. Applying the map~$\pi$, we get $$T_{r_1}T_{r_2} = q^{\ell(r_1)+\ell(r_2)-\ell(r_3)} T_{r_3}.$$
Therefore, it follows from Lemma~\ref{not:fintredmon} that 
$$\begin{array}{ll}T_sT_r = T_{sr},&\textrm{if }s\in S \textrm{ and } \ell(r) =  \ell(r)+1;\\
T_sT_r = qT_{r},&\textrm{if }s\in S\textrm{ and }\ell(sr) = \ell(r);\\ T_eT_r = q^{\ell(r)-\ell(er)}T_{er},&\textrm{if }e\in \crlo.
\end{array}$$     
Assume $s$ lies in $S$ and $r$ lies in~$R$ such that $\ell(sr) = \ell(r)-1$. Denote by $(w_1,e,w_2)$ the normal decomposition of~$r$. By Lemma~\ref{lem:ewwt1}(i), $\ell(sw_1) = \ell(w_1)-1$ and $\ell(sw_1ew_2) = \ell(sw_1)+\ell(w_2)$. Therefore, $T_sT_{w_1} = qT_{sw_1} + (1-q)T_{w_1}$, by~\cite[Theorem~8.4.6]{Gec}, and  $$T_sT_r = T_sT_{w_1}T_{ew_2} = qT_{sw_1}T_{ew_2} + (1-q)T_{w_1}T_{ew_2} = qT_{sr} + (1-q)T_{r}.$$    
\end{proof}
Now, using Theorem~\ref{Th:genThe}, Theorem~\ref{THintro} is a corollary of Theorem~\ref{THintro2}. More precisely,
gathering Corollary~\ref{cor:preshHR} and Theorem~\ref{THintro2}, we get the following result.
\begin{Cor} Let $M$ be a finite reductive monoid over~$\overline{\mathbb{F}}_q$. Consider Notation~\ref{not:fintredmon}. Then the Iwahori-Hecke algebra~$\H(M,B)$ admits the following $\mathbb{C}$-algebra presentation:   
\begin{center}\label{THconclu}
\begin{tabular}{lll}
(HEC1)&$T^2_s = (q-1)T_1+qT_{s}$,&$s\in S$;\\
(HEC2)&$|T_s,T_t\rangle^m = |T_t,T_s\rangle^m$,&$(\{s,t\},m)\in \mathcal{E}(\Gamma)$;\\
(HEC3)&$T_sT_e = T_eT_s$,&  $e\in \crlo$, $s\in \lambda^\star(e)$;\\
(HEC4)&$T_sT_e = T_eT_s = qT_e$,& $e\in \crlo$, $s\in \lambda_\star(e)$;\\
(HEC5)&$T_eT_wT_f = q^{\ell(w)}T_{e\wedge_wf}$,& $e,f\in \crlo$, $w\in \textrm{Red}(e,f)$. 
\end{tabular} 
\end{center}
\end{Cor}
\ \\\ \\\noindent {\bf Acknowledgements.} The author is in debt with F. Digne, B. Leclerc, for useful discussions and with C. Mokler, M. Putcha and L. Renner for helpful email correspondences.


\end{document}